\documentclass{svproc}
\pdfoutput=1

\usepackage{siunitx} %
\usepackage{color}
\usepackage{xcolor}
\usepackage{afterpage}
\usepackage{wrapfig}
\usepackage{float}
 \usepackage{url}

\usepackage{lipsum}
\usepackage{amsfonts}
\usepackage{graphicx}
\usepackage{epstopdf}
\usepackage{tabularx}
\usepackage{mathtools} %
\usepackage{color} %
\usepackage{dashrule} %
\usepackage{enumitem} %
\usepackage{algorithmic}
\ifpdf
  \DeclareGraphicsExtensions{.eps,.pdf,.png,.jpg}
\else
  \DeclareGraphicsExtensions{.eps}
\fi

\definecolor{plot_blue}{RGB}{53, 146, 177}
\definecolor{plot_yellow}{RGB}{253, 209, 7}
\definecolor{plot_red}{RGB}{214, 7, 37}

\title{IMEX Runge-Kutta Parareal for Non-Diffusive Equations}
\author{Tommaso Buvoli \inst{1} and Michael Minion \inst{2}}

\institute{University of California, Merced, Merced CA 95343, USA,\\
\email{tbuvoli@ucmerced.edu},\\
\and
Lawrence Berkeley National Lab,
Berkeley, CA 94720, USA, \\
\email{mlminion@lbl.gov}.\\
}

\usepackage{amsopn}

\ifpdf

\begin{document}

\maketitle

\begin{abstract}
	
	Parareal is a widely studied parallel-in-time method that can achieve meaningful speedup on certain problems. However, it is well known that the method typically performs poorly on non-diffusive equations. This paper analyzes linear stability and convergence for IMEX Runge-Kutta Parareal methods on non-diffusive equations. By combining standard linear stability analysis with a simple convergence analysis, we find that certain Parareal configurations can achieve parallel speedup on non-diffusive equations. These stable configurations all posses low iteration counts, large block sizes, and a large number of processors.  Numerical examples using the nonlinear Schr\"odinger equation demonstrate the analytical conclusions.
\end{abstract}


\def\Oof {{\mathcal{O}}}
\def\Fprop {{\mathcal{F}}}
\def\Gprop {{\mathcal{G}}}
\def\y {{y}}
\def\yn {\y_n}
\def\ynk {{\y_n^k}}
\def\ynkp {{\y_n^{k+1}}}
\def\ynp {\y_{n+1}}
\def\tn  {{t_n}}
\def\tnp  {{t_{n+1}}}
\def\Nstepstot {{N_s}}
\def\Nsteps {{N_T}}
\def\Nprocs {{N_p}}
\def\Nblocks {{N_b}}
\def\NFprop {{N_f}}
\def\NGprop {{N_g}}
\def\Tzero {{T_0}}
\def\Tfinal {{T_{fin}}}
\def\Tblock {{\Delta T}}
\def\Tproc {{\Delta T_p}}
\def\CostG  {{C_\Gprop}}  %
\def\CostF  {{C_\Fprop}}  %
\def\CostB  {{C_B}}  %
\def\CostS  {{C_s}}  %
\def\CostP  {{C_p}}  %
\def\CostGstep  {{c_g}}  %
\def\CostFstep  {{c_f}}  %

\def\Speedup  {{S}}
\def\Eff  {{E}}
\def\u  {{u}}
\def\g  {{u}}
\def\uhat  {{\hat{\u}}}
\def\ghat  {{\hat{\g}}}
\def\uhatk  {{\hat{\u}_k}}
\def\ghatk  {{\hat{\g}_k}}
\def\Tfin  {{T}}
\def\cRate {\|\mathbf{E}\|}
\def\NIters {{K}}
\def\Dt  {{\Delta t}}
\def\DT  {{\Delta T}}
\def\DtF  {{\Delta t_\Fprop}}
\def\DtG  {{\Delta t_\Gprop}}
\def\K  {{K}}
\def\barK  {{\bar{\K}}}
\def\A {\mathcal{A}} %
\section{Introduction}

The numerical solution of ordinary and partial differential equations (ODEs and PDEs) is one of the fundamental tools for simulating engineering and physical systems whose dynamics are governed by differential equations.  Examples of fields where PDEs are used span the sciences from astronomy, biology, and chemistry to zoology, and the literature on methods for ODEs  is well-established (see e.g. 
\cite{Hairer1987-gg,Hairer1991-fj}).  

Implicit-explicit (IMEX) methods are a specialized class of ODE methods that are appropriate for problems where the right-hand side of the equation is additively split into two parts so that %
\begin{align}
    y'(t) = F^E(y,t)+F^I(y,t).
    \label{eq:model_eq}
\end{align}
The key characteristic of an IMEX method is that the term $F^{E}$ (assumed to be non-stiff) is treated explicitly while the term $F^{I}$ (assumed to be stiff) is treated implicitly. In practice, IMEX methods are often used to solve equations that can be naturally partitioned into a stiff linear component that is treated implicitly and a non-stiff nonlinearity that is treated explicitly. %
The canonical example is a nonlinear advection-diffusion type equation, where the stiffness comes from the (linear) diffusion terms while the nonlinear terms are not stiff.  IMEX methods are hence popular in many fluid dynamics settings.  Second-order methods based on a Crank-Nicolson treatment of the diffusive terms and an explicit treatment of the nonlinear terms are a notable example.  However in this study, we consider a different class of problems where the IMEX schemes are applied to a dispersive rather than diffusive term.  Here a canonical example is the non-linear Schr\"odinger equation. Within the class of IMEX methods, we restrict the study here to those based on additive or IMEX Runge-Kutta methods (see e.g., \cite{Ascher1997-er,Calvo2001-ko,Kennedy2003-jd,Boscarino2009-qb}).
In particular, we will study the behavior of the parallel-in-time method, Parareal, constructed from IMEX Runge-Kutta (here after IMEX-RK) methods
applied to non-diffusive problems.

Parallel-in-time methods date back at least to the work of Nievergelt in 1964~\cite{Nievergelt1964} and have seen a resurgence of interest in the last two decades \cite{Gander2015_Review}.  The Parareal method introduced in 2001 \cite{LionsEtAl2001} is perhaps the most well-known parallel-in-time method and can be arguably attributed to catalyzing the recent renewed interest in temporal parallelization.  The emergence of Parareal also roughly coincides with the end of the exponential increase in individual processor speeds in massively parallel computers, a development that has resulted in a heightened awareness of the bottleneck to reducing run time for large scale PDE simulations through spatial parallelization techniques alone.  Although Parareal is a relatively simple method to implement (see Section \ref{sec:parareal}) and can, in principle, be employed using any single-step serial temporal method, one main theme of this paper is that the choice of method is critical to the performance of the algorithm.

Parareal employs a concurrent iteration over multiple time steps to achieve parallel speed up. One of its main drawbacks is that the parallel efficiency is typically modest and is formally bounded by the inverse of the number of iterations required to converge to the serial solution within a pre-specified tolerance. Another well-known limitation is that the convergence of the method is significantly better for purely diffusive problems than for advective or dispersive ones. As we will show, the convergence properties of the Parareal methods considered here are quite complex, and the efficiency is sensitive to the problem being solved, the desired accuracy, and the choice of parameters that determine the Parareal method.  In practice, this makes the parallel performance of Parareal difficult to summarize succinctly. 

Incorporating IMEX integrators into Parareal, enables the creation of new Parareal configurations that have similar stability and improved efficiency compared to Parareal configurations that use fully implicit solvers. IMEX Parareal integrators were first proposed by Wang et al. \cite{Wang2015}, where their stability is studied for equations with a stiff dissipative linear operator and a non-stiff, non-diffusive nonlinear operator. In this work we focus exclusively on non-diffusive equations where the spectrums of $F^E$ and $F^I$ are both purely imaginary. Moreover, we only consider Parareal methods on bounded intervals with a fixed number of iterations. Under these restrictions one can interpret Parareal as a one-step Runge-Kutta method with a large number of parallel stages. By taking this point of view, we can combine classical linear stability and accuracy theory with more recent convergence analyses of Parareal \cite{Ruprecht2018}. Furthermore, fixing the parameters means that the parallel cost of Parareal is essentially known {\it apriori} making comparisons in terms of accuracy versus wall-clock more straightforward.

The main contribution of this work is to introduce new diagrams that combine convergence regions and classical linear stability regions for Parareal on the partitioned Dahlquist problem. The diagrams and underlying analysis can be used to determine whether a particular combination of integrators and parameters will lead to a stable and efficient Parareal method. They also allow us to identify the key Parareal parameter choices that can provide some speedup for non-diffusive problems. Overall, the results can be quite surprising, including the fact that convergence regions do not always overlap with stability regions; this means that a rapidly convergent Parareal iteration does not imply that Parareal  considered as one-step method with a fixed number of iterations is stable in the classical sense.   %

The rest of this paper is organized as follows. In the next section, we present a general overview of IMEX-RK methods and the specific methods used in our study are identified. In Section \ref{sec:parareal}, we provide a short review of the Parareal method followed by a discussion of the theoretical speedup and efficiency.  In Section \ref{sec:stability_convergence}, we conduct a detailed examination of the stability and convergence properties of IMEX-RK Parareal methods. Then, in Section \ref{sec:numerical_experiments} we present several numerical results using the nonlinear Schr\"odinger equation to confirm the insights from the linear analysis. Finally, we present a summary of the findings, along with our conclusions in Section \ref{sec:summary}.
\section{IMEX  Runge-Kutta Methods} \label{sec:imex_rk}

In this section we briefly discuss the IMEX Runge-Kutta (IMEX-RK) methods that are used in this paper. Consider the ODE (\ref{eq:model_eq}) where $F^E$, is assumed to be non-stiff
and while $F^{I}$ is assumed to be stiff. Denoting  $y_{n}$ as the approximation to $y(t_n)$ with $\Dt=\tnp-\tn$,
the simplest IMEX-RK method is forward/backward Euler
  \begin{align}
    y_{n+1}=y_{n}+\Dt \left(F^E(y_{n},t_n) + F^I(y_{n+1},t_{n+1})\right).
  \end{align}
In each step, one needs to evaluate $F^E(y_{n},t_n)$ and then solve  the implicit equation
  \begin{align}
    y_{n+1}-\Dt F^I(y_{n+1},t_{n+1}) =y_{n}+\Dt F^E(y_{n},t_n).
  \end{align}

  Higher-order IMEX methods can be constructed using different families of integrators and  IMEX-RK methods (also called {\it additive} or {\it partitioned}) are one popular choice
(see e.g. \cite{Ascher1997-er,Calvo2001-ko,Kennedy2003-jd,Boscarino2009-qb}).
The generic form for an $s$ stage IMEX-RK method is  
\begin{align}
  y_{n+1} =y_{n}+\Dt \left(\sum_{j=1}^{s} b_j^E F^E(Y_{j},t_n + \Delta t c^{E}_j) + b_j^I F^I(Y_{j},t_n + \Delta t c^{I}_{j})\right),
\end{align}
where the stage values are  
\begin{align}
  Y_{j}=y_{n}+\Dt \left( \sum_{k=1}^{s-1} a_{j,k}^E F^E(Y_k,t_n + \Delta t c^E_k) + \sum_{k=1}^{s} a_{j,k}^IF^I(Y_k, t_n + \Delta t c^I_k) \right).
\end{align}
Such methods are typically encoded using two Butcher tableaus that respectively contain the coefficients $a_{j,k}^E$, $b_{j}^E$, $c_j^{E}$ and $a_{j,k}^I$, $b_{j}^I$, $c_j^{I}$.  As with the Euler method,
each stage of an IMEX method requires the evaluation of $F^E(y_j,t_j)$, $F^I(y_j,t_j)$ and  the solution of the implicit equation
\begin{equation} \label{eq:IMEXstage}
  Y_j- (\Dt a^I_{j,j}) F^I(Y_j,t_n + \Delta t c^I_j) = r_j,         
\end{equation}
where $r_j$ is a vector containing all the known quantities that determine the $j$th stage. IMEX methods are particularly attractive when $F^I(y,t)=Ly$, where $L$ is a linear operator so that (\ref{eq:IMEXstage}) becomes
\begin{equation} \label{eq:IMEXstageL}
  (I- \Dt a^I_{j,j}L)Y_j = r_j.          
\end{equation}
If a fast preconditioner is available for inverting these systems, or if the structure of $L$ is simple, then IMEX methods can provide significant computational savings compared to fully implicit methods.

To achieve a certain order of accuracy, the coefficients $a^E$ and $a^I$ must satisfy both order and matching conditions. Unfortunately, the total number of conditions grows extremely fast with the order of the method, rendering classical order-based constructions difficult. To the best of the authors' knowledge there are currently no IMEX methods beyond order five that have been derived using classical order conditions. However, by utilizing different approaches, such as extrapolation methods \cite{Cardone2013-yg} or spectral deferred correction \cite{Dutt2000SDC,Minion2003IMEX}, it is possible to construct high-order IMEX methods.

In this work, we consider IMEX-RK methods of order one through four. The first and second order methods are the (1,1,1) and (2,3,2) methods from \cite{Ascher1997-er} whose tableaus can be found in Section 2.1 and Section 2.5 respectively. The  third and fourth order methods are the ARK3(2)4L[2]SA and ARK4(3)6L[2]SA from \cite{Kennedy2003-jd}. All the schemes we consider have an L-Stable implicit integrator.

\section{The Parareal method}\label{sec:parareal}

The Parareal method, first introduced in 2001 \cite{LionsEtAl2001}, is a popular approach to time parallelization of ODEs.
In this section, we will give a brief over-view of Parareal and then
present a theoretical model for the parallel efficiency and speedup of the method.

\subsection{Method definition}
\begin{center}
\begin{table}[h]
\begin{tabular}{| c |  l | l |}
 \hline
 Variable & Meaning & Definition  \\ 
 \hline
$\Tfinal$  & Final time of ODE &  Problem specified     \\
$\Nprocs$  & Number of processors & User defined    \\
$\Nstepstot$ & Total Number of fine steps & User defined \\  
$\Nblocks$  & Number of Parareal blocks &  User defined  \\
$K$ & Number of Parareal iterations & User defined or adaptively controlled \\
$\Dt$       & Time step for serial method &  $\Tfinal/\Nstepstot$ \\  
$\Gprop$ & Fine propagator & Here 1 step of IMEX RK method\\  
$\Fprop$ & Fine propagator & $\NFprop$ steps of IMEX RK method\\  
$\NFprop$ & Number of RK steps in $\Fprop$ & $  \Nstepstot/(\Nprocs\Nblocks)$ \\  
$\NGprop$ & Number of RK steps in $\Gprop$ & 1 for all examples\\  
$\Nsteps$ & Total number of fine steps per block &  $\Nstepstot/\Nblocks$ \\  
$\CostS$  & Cost of full serial run  & $\Nstepstot\CostFstep$   \\
$\CostGstep$  & Cost of method per step in $\Gprop$   & User defined   \\
$\CostFstep$  & Cost of method per step in $\Fprop$   & User defined   \\
$\CostF$  & Cost of $\Fprop$   & $\NFprop \CostFstep$   \\
$\CostG$  & Cost of $\Gprop$   & $\NGprop \CostGstep$  \\
$\alpha$  & Ratio of $\Gprop$ to $\Fprop$ cost  & $\CostG/\CostF$   \\
 \hline
\end{tabular}
\caption{Definitions of variable names used in the description of Parareal.}
\end{table}
\end{center}
In the original form, Parareal is straight-forward to describe by a simple iteration.
Let $[\Tzero,\Tfinal]$ be the time interval of interest,
and $\tn$ denote a series of time steps in this interval.  
Next, define coarse and fine propagators $\Gprop$ and $\Fprop$, each of which produces an approximation to the ODE at
$\tnp$ given an approximation to the solution at $\tn$.  

Assume that one has a provisional guess of the solution at each
$\tn$, denoted  $\yn^{0}$. This is usually provided
by a serial application of the coarse propagator $\Gprop$. 
Then the $k$th Parareal iteration is given by
\begin{equation} \label{eq:pararealK}
  \ynp^{k+1} =  \Fprop(\ynk) +  \Gprop(\ynkp) - \Gprop(\ynk),  
\end{equation}
where the critical observation is that the $\Fprop(\yn^k)$ terms can be computed on each time interval in parallel.  
The goal of Parareal is to iteratively compute an approximation to the numerical solution that would result from applying
$\Fprop$ sequentially on $\Nprocs$ time intervals, 
\begin{equation} \label{eq:serialF}
              \ynp =  \Fprop(\yn),  \hbox{ \hspace{7mm} for \hspace{1mm} }  n=0 \ldots \Nprocs-1
\end{equation}
where each interval is assigned to a different processor.
As shown below, assuming that
$\Gprop$ is computationally much less expensive than $\Fprop$ and that the method converges in few enough iterations, parallel speedup can be obtained. 

Part of the appeal of the Parareal method is that the propagators $\Gprop$ and $\Fprop$ are not constrained by the definition of the method.  Hence, parareal as written can in theory be easily implemented using any numerical ODE method for $\Gprop$ and $\Fprop$.  Unfortunately, as discussed below, not all choices lead to efficient or even convergent parallel numerical methods, and the efficiency of the method is sensitive to the choice of parameters.

Note that as described, the entire Parareal method can be considered
as a self-starting, single step method for the interval
$[\Tzero,\Tfinal]$ with time step $\DT = \Tfinal-\Tzero$.  In the
following Section, the classical linear stability of Parareal as a
single step method will be considered for $\Gprop$ and $\Fprop$ based
on IMEX-RK integrators.  This perspective also highlights the
fact that there is a choice that must be made for any particular
Parareal run regarding the choice of $\DT$.  To give a concrete
example for clarity, suppose the user has an application requiring
1024 time steps of  some numerical method to compute the desired
solution on the time interval $[0,1]$, and that 8 parallel processors
are available.  She could then run the Parareal algorithm on 8
processors with 128 steps of the serial method corresponding to
$\Fprop$.  Alternatively, Parareal could be run as a single step
method on two {\it blocks} of time steps corresponding to $[0,1/2]$ and $[1/2,1]$
with each block consisting
to 512 serial fine time steps, or 64 serial steps corresponding to $\Fprop$ for each processor
on each block.  These two blocks would necessarily be computed
serially with the solution from the first block at $t=1/2$ serving as the initial
initial condition on the second block. %

\subsection{Cost and theoretical parallel speedup and efficiency} 

We describe a general framework for estimating the potential speedup
for the Parareal method in terms of the reduction in the run time of the method.  Although theoretical cost estimates have been considered in detail before (see e.g. \cite{Aubanel2011-rc}, we repeat the basic derivation for the specific assumptions of the IMEX-RK based methods used here and for the lesser known estimates for multiple block Parareal methods. 
For simplicity, assume that an initial value
ODE is to be solved with some method 
requiring $\Nstepstot$ time steps to complete the simulation on the interval
$[0,\Tfinal]$.  We assume further that the same method will be used
in the fine propagator in Parareal.
If each step of the serial method has cost $\CostFstep$, then the total serial cost is
\begin{equation}
  \label{eq:serial_cost}
  \CostS= \Nstepstot \CostFstep.
\end{equation}
In the numerical examples present in section 5, both the coarse and fine propagators consist
of a number of steps of an IMEX-RK method applied to a
pseudo-spectral discretization of a PDE. The main cost in general is
then the cost of the FFT used to compute explicit nonlinear spatial function evaluations.
Hence each step of either IMEX-RK method has essentially a
fixed cost, denoted $\CostFstep$ and $\CostGstep$.  This is in contrast to the case where implicit equations
are solved with an iterative method and the cost per time step
could vary  considerably by step.

Given $\Nprocs$ available processors, the Parareal
algorithm can be applied to $\Nblocks$ blocks of time intervals, with
each block having length $\Tblock=\Tfinal/\Nblocks$.  Again for
simplicity we assume that in each time block, each processor is
assigned a time interval of equal size $\Tproc=\Tblock/\Nprocs$.
Under these assumptions, $\Fprop$ is now determined to be
$\NFprop=\Nstepstot/(\Nprocs \Nblocks)$ steps of the serial method.
Parareal is then defined by the choice of 
$\Gprop$, which we assume here is constant across processors and
blocks consisting of
$\NGprop$ steps of either the same, or different RK method as used in
$\Fprop$ with cost per step $\CostGstep$.
Let  $\CostF=\NFprop \CostFstep $ be  the time needed to compute $\Fprop$,
and likewise, let $\CostG = \NGprop \CostGstep$ be the cost of the coarse propagator.

The cost of $K$ iterations of Parareal performed on
a block is the sum of the cost of the predictor on a block, 
$\Nprocs  \CostG$, plus the additional cost of each iteration.
In an ideal setting where each processor computes a quantity as soon
as possible and communication cost is neglected, the latter is simply
the $K(\CostF + \CostG)$.  Hence the total cost of Parareal on a block is
\begin{equation}
  \label{eq:parareal_cost2}
  \CostB=  \Nprocs \CostG+ K(\CostF + \CostG).
\end{equation}
The total cost of Parareal is the sum over blocks
\begin{equation}
  \CostP =\sum_{i=1}^\Nblocks\left( \Nprocs \CostG+ K_i(\CostF + \CostG)\right) = \Nblocks\Nprocs\CostG +  (\CostF + \CostG) \sum_{i=1}^\Nblocks  K_i,
\end{equation}
where $K_i$ is the number of iterations required to converge on block $i$.
Let $\barK$ denote the average number of iterations across the blocks
then
\begin{equation}
  \CostP = \Nblocks \left(\Nprocs\CostG +  \barK (\CostF + \CostG)\right).
\end{equation}
Note the first term 
$\Nblocks \Nprocs\CostG$ is exactly the 
cost of applying the coarse propagator over the entire time interval.

Finally, denoting $\alpha=\CostG/\CostF$,  the speedup $\Speedup=\CostS/\CostP$ is then
\begin{equation}
  \label{eq:parareal_speedup2}
  \Speedup= \frac{\Nprocs}{\Nprocs \alpha + \barK (1+ \alpha )}.
\end{equation}
For fixed $K$, where $\barK=K$, this reduces to the usual estimate (see e.g.Eq. (28) \cite{Minion2010} or Eq. (19) of \cite{Aubanel2011-rc}).
The parallel efficiency of parareal $\Eff=\Speedup/\Nprocs$ is
\begin{equation}
  \label{eq:parareal_efficiency2}
  \Eff= \frac{1}{(\Nprocs+\barK)\alpha + \barK}.
\end{equation}

A few immediate observations can be made from the formulas for
$\Speedup$ and $\Eff$. Clearly the bound on efficiency is  $\Eff < 1/\barK$.
Further, if significant speedup is to be achieved, it should be true that  $\barK$ is significantly less than $\Nprocs$ and
$\Nprocs \alpha$ is small as well.
As will be demonstrated later, the total number of Parareal iterations required is
certainly problem dependent and also dependent on the choices of
$\Fprop$ and $\Gprop$.
It might seem strange at first glance that the number of blocks chosen does not appear explicitly in the above formulas for
$\Speedup$ and $\Eff$.  Hence it would seem better to choose more blocks of shorter length so that $\barK$ is minimized.
Note however that increasing the number of blocks by a
certain factor with the number of processors fixed means that $\NFprop$ will decrease by the same factor.  If the cost of
the coarse propagator $\CostG$ is independent of the number of blocks (as in the common choice of $\Gprop$ being a single
step of a given method, i.e., $\NGprop=1$), then $\alpha$ will hence increase by the same factor.
Lastly,  one can derive the total speedup by also considering the speedup over each block, $\Speedup_i$ as
\begin{align}
     \Speedup=\frac{1}{\sum_{i=1}^\Nblocks \frac{1}{\Nblocks \Speedup_i}}=\frac{1}{\frac{1}{\Nblocks}\sum_{i=1}^\Nblocks \frac{1}{\Speedup_i}}.
 \end{align}

Finally, we should note that more elaborate parallelization strategies than that discussed above are possible, for example \cite{Aubanel2011-rc,BerryEtAl2012,ArteagaEtAl2015,Ruprecht2017_lncs}.
\definecolor{plot_grey}{RGB}{80 80 80}
\definecolor{plot_blue_light}{RGB}{84 127 255}
\definecolor{plot_blue_dark}{RGB}{57 80 151}
\definecolor{plot_yellow}{RGB}{255 186 65}

\section{Non-diffusive Dalquist: stability, convergence, accuracy }
\label{sec:stability_convergence}

In this section, we analyze linear stability and convergence  properties for IMEX-RK Parareal methods for non-diffusive problems.
There have been multiple previous works that have analyzed  convergence and stability properties of Parareal. Bal \cite{Bal2005} analyzed Parareal methods with fixed parameters, and Gander and Vandewalle \cite{GanderVandewalle2007_SISC} studied the convergence of parareal on both bounded and unbounded intervals as the iterations $k$ tends to infinity. More recently, Southworth et al. \cite{Southworth2019,southworth2020tight} obtained tight convergence bounds for Parareal applied to linear problems. Specific to stability for non-diffusive equations, Staff and Ronquist \cite{StaffRonquist2005} conducted an initial numerical study, Gander \cite{Gander2008} analyzed the stability of Parareal using characteristics, and Ruprecht \cite{Ruprecht2018} studied the cause of instabilities for wave equations. In this work we will use the work of Ruprecht as a starting point to study stability and convergence for Parareal integrators that are constructed using IMEX-RK integrators.

In this work, we consider the non-diffusive partitioned Dahlquist test problem  
\begin{align}
	\left\{
	\begin{aligned}
		& y' = i \lambda_1 y + i\lambda_2y \\
		& y(0) = 1	
	\end{aligned}
	\right.
	&& \lambda_1, \lambda_2 \in \mathbb{R},
	\label{eq:dispersive-dahlquist}
\end{align} 
where the term $i\lambda_1 y$ is treated implicitly and the term $i\lambda_2y$ is treated explicitly. This equation is a generalization of the Dahlquist test problem that forms the basis of classical linear stability theory \cite[IV.2]{wanner1996solving}, and the more general equation with $\lambda_1, \lambda_2 \in \mathbb{C}$ has been used to study the stability properties of various specialized integrators \cite{ascher1995implicit,izzo2017highly,sandu2015generalized,cox2002ETDRK4,krogstad2005IF,buvoli2020class}. In short, (\ref{eq:dispersive-dahlquist}) highlights stability for (\ref{eq:model_eq}) when $F^{E}(y,t)$ and $F^{I}(y,t)$ are autonomous, diagonalizable linear operators that share the same eigenvectors, and have a purely imaginary spectrum.

When solving (\ref{eq:dispersive-dahlquist}), a classical one-step integrator (e.g. an IMEX-RK method) reduces to an iteration of the form
\begin{align}
	y_{n+1} = R(iz_1,iz_2) y_n	 \quad \text{where} \quad z_1 = h\lambda_1,~ z_2 = h \lambda_2,
	\label{eq:one-step-dahlquist-iteration}
\end{align}
and $R(\zeta_1,\zeta_2)$ is the stability function of the method. A Parareal algorithm over an entire block can also be interpreted as a one-step method that advances the solution by $\Nprocs$ total time steps of the integrator $\mathcal{F}$. Therefore, when solving (\ref{eq:dispersive-dahlquist}), it reduces to an iteration of the form
\begin{align}
	y_{(\Nprocs(n+1))} = R(iz_1,iz_2) y_{(\Nprocs n)}.
	\label{eq:one-step-dahlquist-iteration-parareal}
\end{align}
The stability function $R(\zeta_1, \zeta_2)$ plays an important role for both convergence and stability of Parareal, and the approach we take for determining the stability functions and convergence rate is identical to the one presented in \cite{Ruprecht2018}. 

The formulas and analysis presented in the following two subsections pertain to a single Parareal block. Since we will compare Parareal configurations that vary the number of fine steps $\NFprop$ (so that the fine integrator $\mathcal{F}$ is $\mathcal{F} = f^\NFprop$), it is useful to introduce the {\em blocksize} $\Nsteps = \Nprocs \NFprop$ which corresponds to the total number of steps that the integrator $f$ takes over the entire block.%

\subsection{Linear stability}

The stability region for a one-step IMEX method with stability function $R(\zeta_1, \zeta_2)$ is the region of the complex $\zeta_1$ and $\zeta_2$ plane given by
\begin{align}
	\hat{\mathcal{S}} = \left\{ (\zeta_1, \zeta_2) \in \mathbb{C}^2 : |R(\zeta_1, \zeta_2)| \le 1 \right\}.	
\end{align}
Inside $\hat{\mathcal{S}}$ the amplification factor $|R(\zeta_1, \zeta_2)|$ is smaller than or equal to on, which ensures that the timestep iteration remains forever bounded. For traditional integrators one normally expects to take a large number of timesteps, so even a mild instability will eventually lead to unusable outputs. 

The full stability region $\hat{\mathcal{S}}$ is four dimensional is difficult to visualize. Since we are only considering the non-diffusive Dahlquist equation, we restrict ourselves to the simpler two dimensional stability region
\begin{align}
	\mathcal{S} = \left\{ (z_1, z_2) \in \mathbb{R}^2 : |R(iz_1, iz_2)| \le 1 \right\}.	
	\label{eq:linear-stability-region} 
\end{align}
Moreover, all the integrators we consider have stability functions that satisfy 
\begin{align}
	R(iz_1, iz_2) = R(-iz_1, -iz_2)
	\label{eq:stability-symmetry-condition}
\end{align}
which means that we can obtain all the relevant information about stability by only considering $\mathcal{S}$ for $z_1 \ge 0$.
 
 \subsubsection{Linear stability for IMEX-RK}

\begin{figure}[h!]

	\includegraphics[width=\linewidth,trim={2cm 1.5cm 2cm 1.5cm},clip]{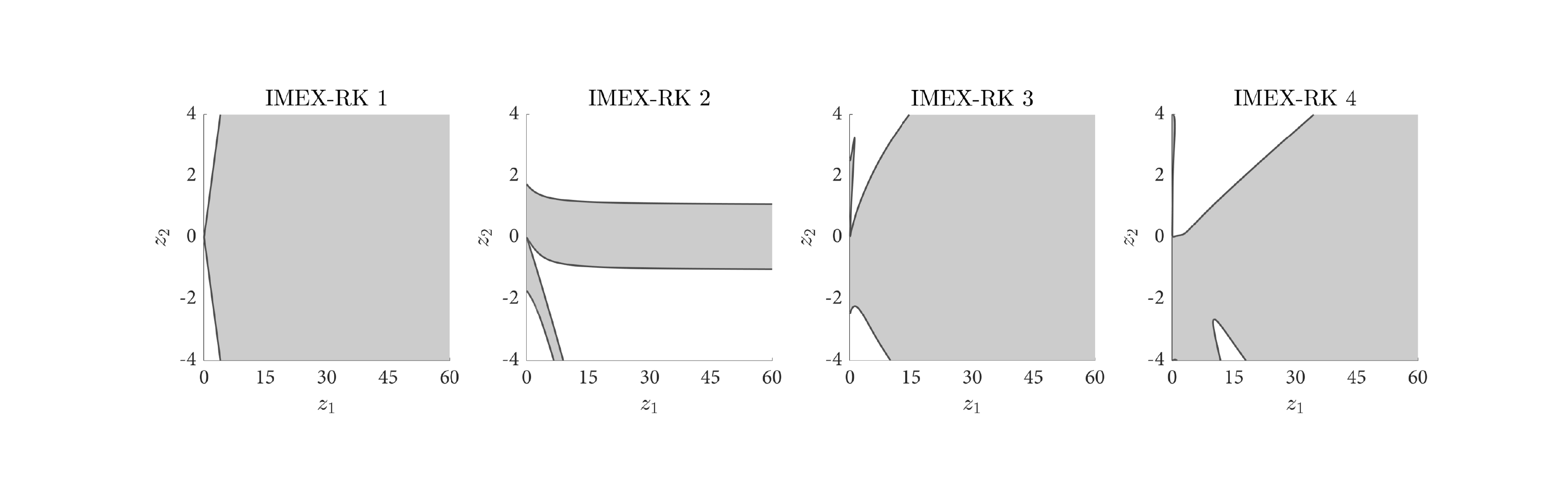}
	
	\vspace{0.5em}
	\includegraphics[width=\linewidth,trim={2cm 1cm 2cm 0cm},clip]{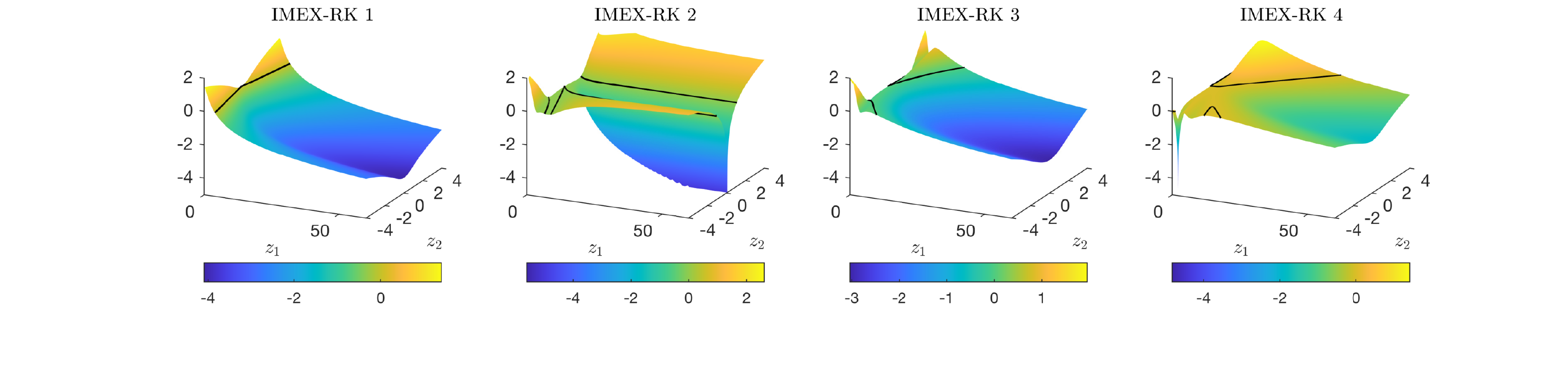}
	
	\caption{non-diffusive linear stability regions (\ref{eq:linear-stability-region}) for IMEX-RK methods (top) and surface plots showing $\log(\text{amp}(iz_1,iz_2))$ (bottom). For improved readability we scale the $z_1$ and $z_2$ axes differently. For the amplitude function plots, zero marks the cutoff for stability since we are plotting the log of the amplitude function.
	}
	\label{fig:imrk-stability}

\end{figure}

Before introducing stability for Parareal, we briefly discuss the linear stability properties of the four IMEX-RK methods considered in this work. In Figure \ref{fig:imrk-stability} we present two-dimensional stability regions (\ref{eq:linear-stability-region}) and surface plots that show the corresponding amplitude factor.  When $z_2=0$, IMEX-RK integrators revert to the fully implicit integrator. Since the methods we consider are all constructed using an  L-stable  implicit method, the amplification factor will approach to zero as $z_1 \to \infty$. This implies that we should not expect good accuracy for large $|z_1|$ since the exact solution of the non-diffusive Dahlquist equation always has magnitude one. As expected, this damping occurs at a slower rate for the more accurate high-order methods.

 \subsubsection{Linear stability for Parareal}

 The importance of linear stability for Parareal (i.e the magnitude of $R(z_1, z_2)$ from (\ref{eq:one-step-dahlquist-iteration-parareal})) depends on the way the method is run and on the severity of any instabilities. In particular, we consider two approaches for using Parareal. In the first approach, one fixes the number of processors and integrates in time using multiple Parareal blocks. This turns Parareal into a one-step RK method; therefore, if one expects to integrate over many blocks, then the stability region becomes as important as it is for a traditional integrator. 
 
 An alternative approach is to integrate in time using a single large Parareal block. If more accuracy is required, then one simply increases the number of timesteps and/or processors, and there is never a repeated Parareal iteration. In this second scenario we can relax traditional stability requirements since a mild instability in the resulting one-step Parareal method will still produce usable results. However, we still cannot ignore large instabilities that amplify the solution by multiple orders of magnitude.
 
To analyze the linear stability of parareal, we first require a formula for its stability function. In \cite{Ruprecht2018}, Ruprecht presents a compact formulation for the stability function of a single Parareal block. He first defines the matrices
	\begin{align}
	\mathbf{M}_F = 
		\left[
			\begin{array}{cccc}
			I \\
			-F & I \\
			   & \ddots & \ddots \\
			   & & -F & I
			\end{array}
		\right] &&
	\mathbf{M}_G = 	\left[
			\begin{array}{cccc}
			I \\
			-G & I \\
			   & \ddots & \ddots \\
			   & & -G & I
			\end{array}
		\right],
		\label{eq:parareal-matrices}
	\end{align}
where the constants $F = R_f(iz_1, iz_2)^{\NFprop}$ and $G = R_c(iz_1, iz_2)^{\NGprop}$ are the stability functions for the fine propagator $\Fprop$ and the coarse propagator $\Gprop$. The stability function for Parareal is then
\begin{align}
	R(iz_1, iz_2) = \mathbf{c}_2 \left( \sum_{j=0}^k \mathbf{E}^j \right) \mathbf{M}_G^{-1} \mathbf{c}_1,
\end{align}
where $\mathbf{E} = \mathbf{I} - \mathbf{M}_G^{-1}\mathbf{M}_F^{-1}$ and $\mathbf{c}_1 \in \mathbb{R}^{\Nprocs+1}$, ${\mathbf{c}_2 \in \mathbb{R}^{1,\Nprocs+1}}$ are $\mathbf{c_1} = \left[ 1, 0, \ldots, 0 \right]^T$ and $\mathbf{c_2} = \left[ 0, \ldots, 0, 1 \right]$.

\subsection{Convergence}
\label{subsec:convergence}
 A Parareal method will always converge to the fine solution after $\Nprocs{}$ iterations. However, to obtain parallel speedup, one must achieve convergence in substantially fewer iterations.  Convergence rates for a linear problem can be studied by writing the Parareal iteration in matrix form, and computing the maximal singular values of the iteration matrix \cite{Ruprecht2018}. Below, we summarize the key formulas behind this observation.
 
For the linear problem (\ref{eq:dispersive-dahlquist}) the Parareal iteration (\ref{eq:pararealK}) reduces to 
\begin{align}
\mathbf{M}_G \mathbf{y}^{k+1} = (\mathbf{M}_G - \mathbf{M}_F)\mathbf{y}^{k} + \mathbf{b}
\end{align}
where $\mathbf{y} = [y_0, y_{1}, y_{2}, \ldots, y_{\Nprocs{}}]^T$ is a vector containing the approximate Parareal solutions at each fine timestep of the integrator $\mathcal{F}$, the matrices $\mathbf{M}_G, \mathbf{M}_F \in \mathbb{R}^{\Nprocs{}+1,\Nprocs{}+1}$ are defined in (\ref{eq:parareal-matrices}) and the vector $\mathbf{b} \in \mathbb{R}^{\Nprocs{}+1}$ is $\left[y_0, 0, \ldots, 0\right]^{T}$.
The Parareal algorithm can now be interpreted as a fixed point iteration that converges to the fine solution 
\begin{align}
	\mathbf{y}_F = \left[1, F, F^2, \ldots, F^\Nprocs\right]^T y_0
\end{align}
and whose error $\mathbf{e}_k = \mathbf{y}^k - \mathbf{y}_F$ evolves according to
	\begin{align}
		\mathbf{e}^k = \mathbf{E} \mathbf{e}^{k-1}	\quad \text{where} \quad \mathbf{E} = \mathbf{I} - \mathbf{M}_G^{-1} \mathbf{M}_F.
		\label{eq:iteration_matrix_parareal_error}
	\end{align}
	Since Parareal converges after $\Nprocs{}$ iterations, the matrix $\mathbf{E}$ is nilpotent and convergence rates cannot be understood using the spectrum. However, monotonic convergence is guaranteed if $\|\mathbf{E}\| < 1$ since
	\begin{align*}
			\|e^{k+1}\| \le \| \mathbf{E}\| \|e^k\| <  \|e^k\|,
	\end{align*} 
	where $\| \cdot \|$ represents any valid norm. We therefore introduce the {\em convergence region}
	\begin{align}
	 	\mathcal{C}_p = \left\{ (z_1, z_2): \| E \|_p < 1\right\}	
	 	\label{eq:convergence-region}
	\end{align}
	that contains the set of all $z_1$, $z_2$ where the $p$-norm of $\mathbf{E}$ is smaller than one and the error iteration (\ref{eq:iteration_matrix_parareal_error}) is contractive.	 Note that for rapid convergence that leads to parallel speedup one also needs $\|\mathbf{E}\|_p \ll 1$.

	\subsubsection{Two-norm for bounding $\mathbf{E}$} In \cite{Ruprecht2018}, Ruprecht selects ${\| \mathbf{E}\|_2 = \max_j \sigma_j}$, where $\sigma_j$ is the $j$th singular value of $\mathbf{E}$. However, the two-norm needs to be computed numerically, which prevents us from understanding the conditions that guarantee fast convergence.	 
	 
	 \subsubsection{Infinity-norm for bounding $\mathbf{E}$} 
	
	 If we consider the $\infty$-norm, we can exploit the simple structure of the matrix $\mathbf{E}$ to obtain the exact formula
	\begin{align}
		\| \mathbf{E} \|_\infty = \frac{1 - |G|^{\Nprocs}}{1 - |G|} |G - F|.
		\label{eq:parareal_iteration_inf_norm}
	\end{align}
	This equality can be obtained directly through simple linear algebra (See Appendix \ref{ap:inf_norm_parareal_iteration_matrix}) and is  similar to the formula used in more sophisticated convergence analysis of Parareal \cite{GanderVandewalle2007_SISC,Southworth2019} and MGRIT \cite{DobrevEtAl2017}. By using this exact formula, we can understand the requirements that must be placed on the coarse and fine integrators to guarantee a rapidly convergent Parareal iteration. We summarize them in three remarks.	
	\begin{remark}
		If $G$ is stable, so that its stability function is less then one, and $|G - F| < \frac{1}{\Nprocs}$ then the Parareal iteration converges monotonically. Notice that when $|G|<  1$, then 
		\begin{align}
			\frac{1 - |G|^{\Nprocs}}{1 - |G|} = \sum_{j=0}^{\Nprocs - 1} |G|^j < \Nprocs.
		\end{align}
		Therefore from (\ref{eq:parareal_iteration_inf_norm}) it follows that if $|G - F| < \frac{1}{\Nprocs}$ then $\|\mathbf{E}\|_{\infty} < 1$. Note however, that this is not always mandatory, and in the subsequent remark we show how this restriction is only relevant for modes with no dissipation. Nevertheless, if we want to add many more processors to a Parareal configuration that converges for all modes, then we also require a coarse integrator that more closely approximates the fine integrator. One way to satisfy this restriction is by keeping $\Nsteps$ fixed while increasing the number of processors; this shrinks the stepsize of the coarse integrator so that it more closely approximates the fine integrator. Another option is to simply select a more accurate coarse integrator or increase the number of coarse steps $\NGprop$.
	\end{remark}

	\begin{remark}
		\label{remark:parareal_dispersive_dissipative}
		It is more difficult to achieve large convergence regions for a non-diffusive equation than for a diffusive one. If we are solving a heavily diffusive problem $y'=\rho_1y + \rho_2y$ where $\text{Re}(\rho_1+\rho_2) \ll 0$ with an accurate and stable integrator, then $|G| \ll 1$. Conversely, if we are solving a stiff non-diffusive problem (\ref{eq:dispersive-dahlquist}) with an accurate and stable integrator we expect that $|G| \sim 1$. Therefore
			\begin{align}
				\frac{1 - |G|^{\Nprocs}}{1 - |G|} \sim 
				\begin{cases}
					1 & \text{Diffusive Problem,} \\
					\Nprocs & \text{Non-Diffusive Problem.}
				\end{cases}
			\end{align}
			From this we see that the non-diffusive case is inherently more difficult since we require that the difference between the coarse integrator and the fine integrator should be much smaller than $\frac{1}{\Nprocs}$ for fast convergence. Moreover, any attempts to pair an inaccurate but highly stable coarse solver ($|G| \ll 1$) with an accurate fine solver ($|F| \sim 1$) will at best lead to slow convergence for a non-diffusive problem since $|G - F| \sim 1$. Rapid convergence is possible if both $|F| \ll 1$ and $|G| \ll 1$, however this is not meaningful convergence since both the coarse and fine integrator are solving the non-diffusive problem inaccurately.
	\end{remark}

	\begin{remark}
		If $G$ is not stable (i.e. $|G|>1$), then fast convergence is only possible if $F$ is also unstable so that $|F|>1$. Convergence requires that the difference between the coarse and fine iterator is sufficiently small so that 

			\begin{align}
				|G - F| < \frac{1 - |G|}{1 - |G|^{\Nprocs}}.
			\end{align}
			Since $G$ and $F$ are complex numbers we can interpret $\frac{1 - |G|}{1 - |G|^{\Nprocs}}$ as the radial distance between the numbers. If we want $|F|\le 1$, then $G$ can never be more than the distance $\frac{1 - |G|}{1 - |G|^{\Nprocs}}$ from the unit circle. Therefore we require that
			\begin{align}
				|G| - \frac{1-|G|}{1-|G|^{\Nprocs}} \le 1 \quad \implies \quad  |G|\le 1.
			\end{align}

	\end{remark}

\subsection{Linear stability and convergence plots for Parareal}

The aim of this subsection is to broadly categorize the effect that each of the Parareal parameters have on stability and convergence. The first parameters that we consider are the coarse and fine integrator. Since we are considering IMEX-RK methods with orders one to four, there are ten possible IMEX-RK pairings where the fine integrator has higher or equivalent order compared to the coarse integrator. The remaining parameters are the number of processors, the number of Parareal iterations $\NIters$, and the number of coarse and fine integrator steps, $\NGprop$ and $\NFprop$. 

In practice, Parareal is commonly run with an adaptively selected $\NIters$ that causes the method to iterate until a pre-specified residual tolerance is satisfied. When discussing linear stability we instead assume that Parareal always performs a fixed number of iterations. This simplifies our analysis and  allows us to quantify how  iteration count affects stability. Since the matrix $\mathbf{E}$ from (\ref{eq:iteration_matrix_parareal_error}) does not depend on $\NIters$, the number of iterations has no effect on the convergence regions. Therefore, our convergence results apply to both non-adaptive and adaptive Parareal.

To reduce the number of total parameters further, we always take the number of coarse integrator steps $\NGprop$ to be one. Even with this simplification there are still too many degrees of freedom to discuss all the resulting Parareal methods. Therefore, we only show several example plots that capture the essential phenomena, and provide a set of general remarks to encapsulate our main observations. Those who wish to see additional stability and convergence plots, can download our Matlab code \cite{t_buvoli_2021_4513662} which can be used to generate figures for any other set of Parareal parameters.

To showcase the stability and convergence properties of Parareal, we present two dimensional plots that overlay the linear stability region (\ref{eq:linear-stability-region}) and the convergence region (\ref{eq:convergence-region}) with $p=\infty$. We start with a simple example that demonstrates how we construct our plots. Consider the Parareal integrator with IMEX-RK3 and IMEX-RK4 as the coarse and fine integrator, respectively, and with parameters $\NIters=2$, $\NGprop=1$, $\NFprop=8$, $\Nprocs=64$. The following three figures show the convergence region (left), the stability region (middle) and an overlay of both (right). Each plot shows the ($z_1$, $z_2$) plane where we only consider $z_1 \ge 0$ due to the symmetry condition (\ref{eq:stability-symmetry-condition}). 
	\begin{center}
		\includegraphics[width=0.32\linewidth]{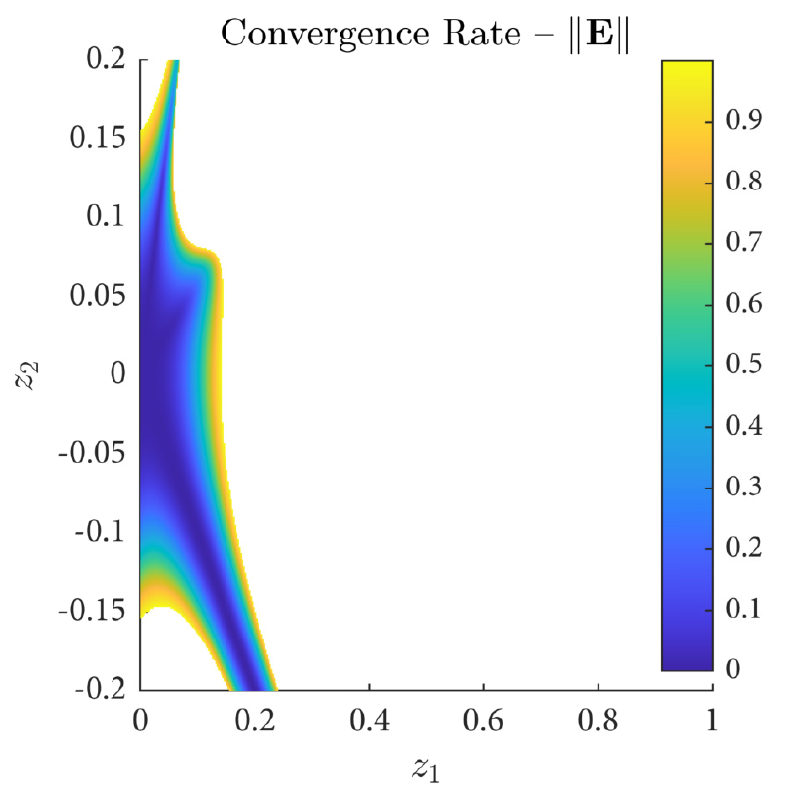}
		\includegraphics[width=0.32\linewidth]{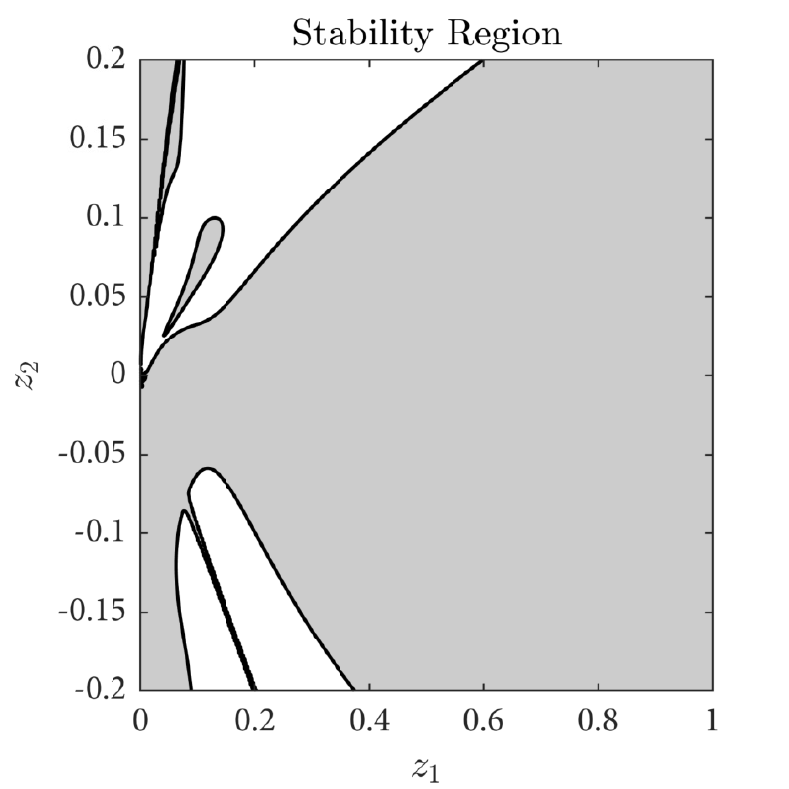}
		\includegraphics[width=0.32\linewidth]{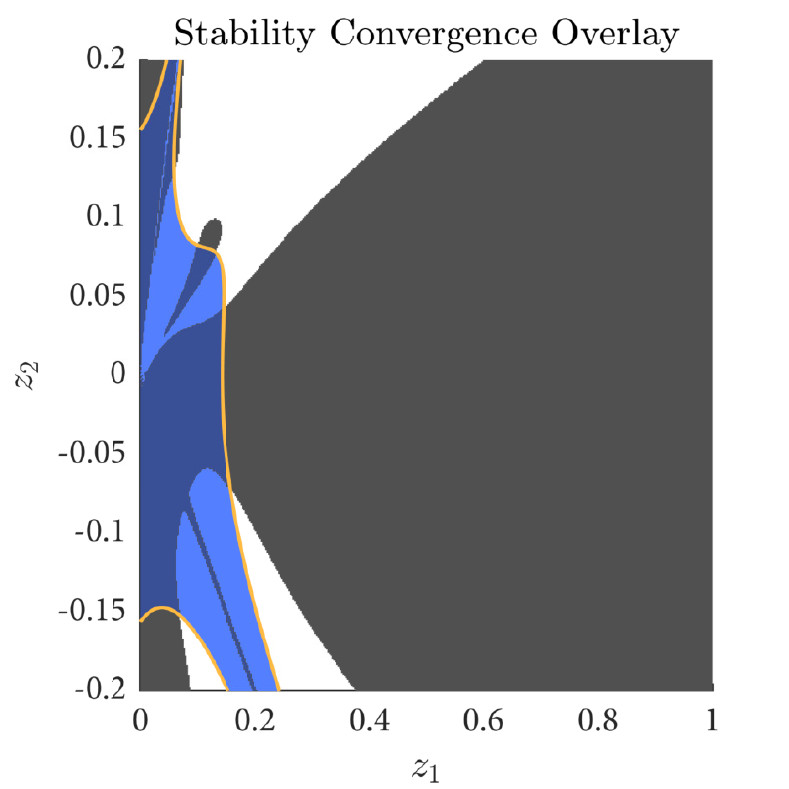}	

			\vspace{0.5em}
			\begin{small}
			\renewcommand*{\arraystretch}{1.25}
			\begin{tabular}{lll}
			{\bf Colo}r & {\bf Set} & {\bf Description} \\ \hline \\[-.9em]
				{\colorbox{plot_blue_light}{\vphantom{X}\hspace{1em}}} Light Blue~~ & $\mathcal{C}_\infty \setminus \mathcal{S}$ & Contractive iteration, unstable one-step method. \\[.2em] 
				{\colorbox{plot_blue_dark}{\vphantom{X}\hspace{1em}}} Dark Blue & $\mathcal{S} \cap \mathcal{C}_\infty$ & Contractive iteration, stable one-step method.\\[.2em] 
				{\colorbox{plot_grey}{\vphantom{X}\hspace{1em}}} Dark Gray & $\mathcal{S} \setminus \mathcal{C}_\infty$ & Non-contractive iteration, stable one-step method. \\[.2em]
				{\colorbox{plot_yellow}{\vphantom{X}\hspace{1em}}} Yellow & $\cRate_{\infty} = 1$~~ & Boundary of convergence region.
			\end{tabular}
			\end{small}
		\end{center}

			\noindent The convergence rate plot has a colorbar that corresponds to the norm of the Parareal iteration matrix $\mathbf{E}$ from (\ref{eq:iteration_matrix_parareal_error}). In the overlay plot there are three distinct regions that are colored according to the legend shown above.
		
	To simplify the comparison between different Parareal methods, we scale all the stability functions used to compute both stability and convergence relative to the number of total fine integrator steps $\Nsteps$; in other words, our plots are generated using the scaled stability function $\widehat{R}(iz_1, iz_2) = R(i\Nsteps z_1, i\Nsteps z_2)$.

\begin{figure}[h!]
	\centering
	\includegraphics[width=0.98\linewidth,trim={10cm 20cm 8cm 9cm},clip]{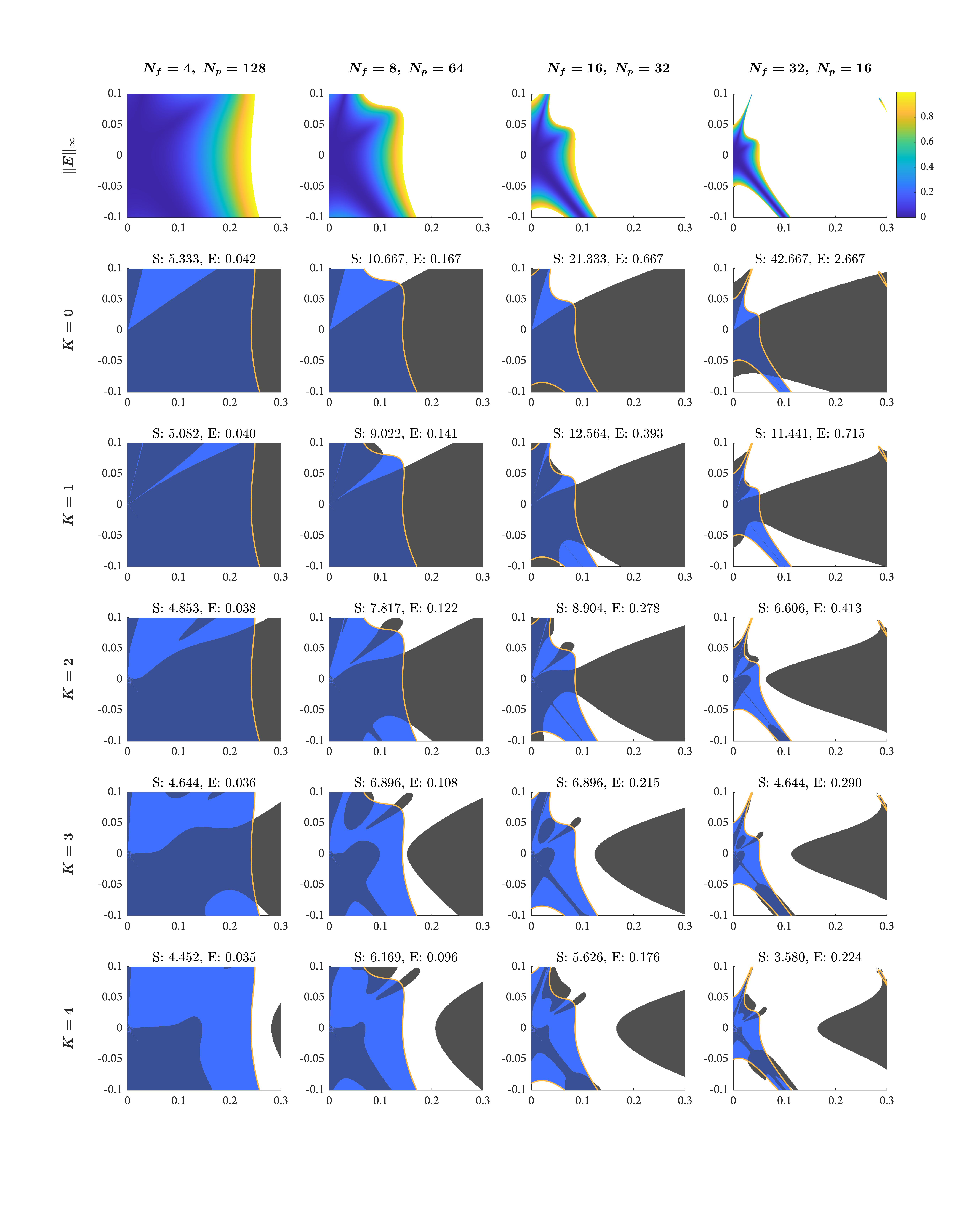}
	
	\caption{Stability and convergence overlay plots for Parareal configurations with a block size of $\Nsteps=512$ and IMEX-RK3, IMEX-RK4 as the coarse and fine integrator. The number of processors increases as one moves rightwards in the horizontal direction. The top row shows the convergence rate plot for the method, and all subsequent rows show stability convergence overlays for an increasing number of iterations. The subfigure titles show the theoretical speedup  (S) and efficiency (E) for each Parareal configuration.}
	\label{fig:imex_parareal_512_close}	
\end{figure}

\begin{figure}[h!]
	\centering
	\includegraphics[width=0.98\linewidth,trim={10cm 20cm 8cm 9cm},clip]{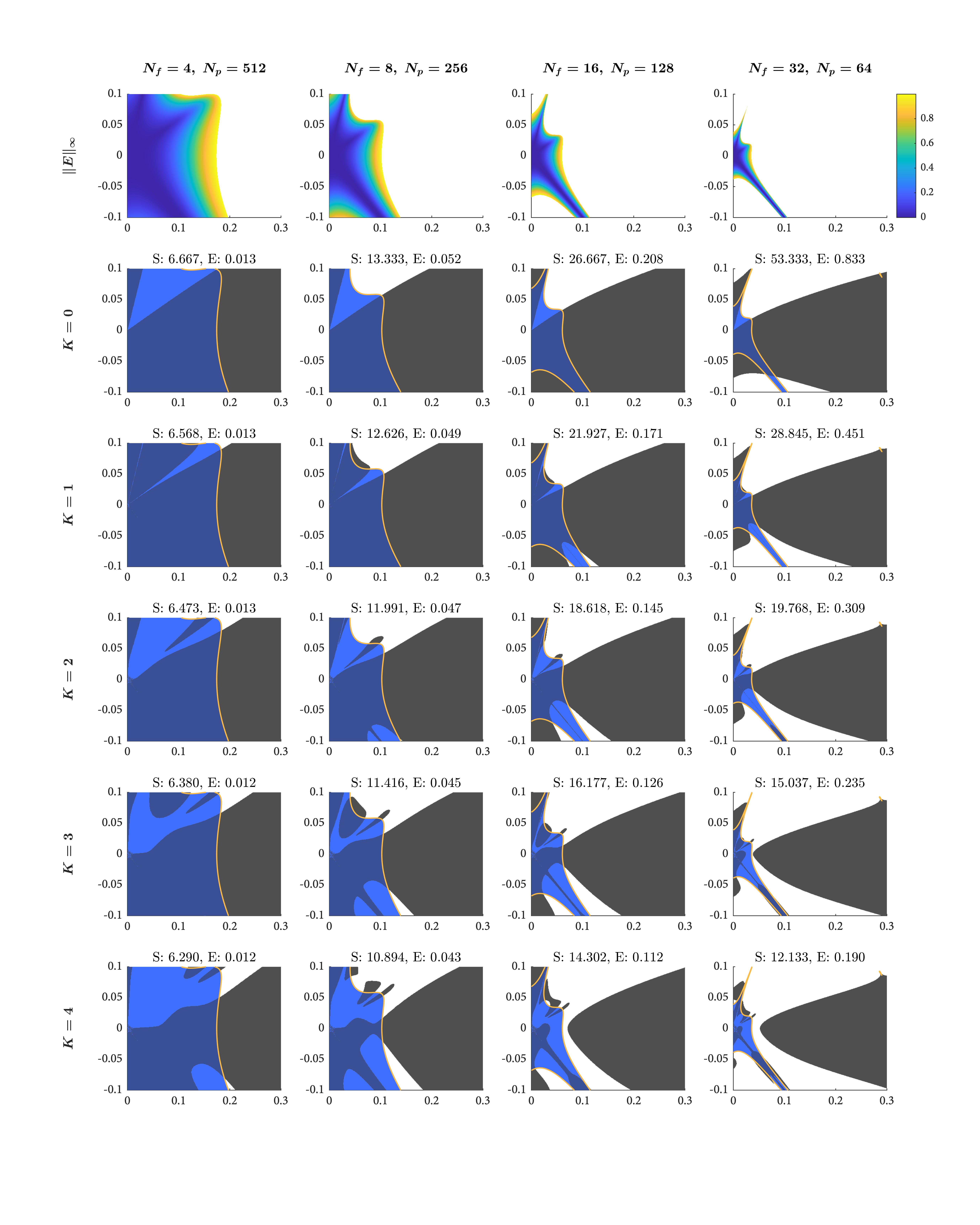}
	
	\caption{Stability and convergence overlay plots for Parareal configurations with a block size of $\Nsteps=2048$ and IMEX-RK3, IMEX-RK4 as the coarse and fine integrator.}
	\label{fig:imex_parareal_2048_close}	
\end{figure}

	For brevity we only present plots for Parareal methods with IMEX-RK3 and IMEX-RK4 as the coarse and fine integrator, and with the following parameters:
	\begin{align*}
		\Nsteps \in \{ 512, 2048 \}, \quad \NGprop \in \{1 \}, \quad \NFprop \in \{4, 8, 16, 32\}, \quad k \in \{1, 2, 3, 4\}.
	\end{align*}
	For each configuration the number of processors $\Nprocs = \Nsteps / \NFprop$. In Figures \ref{fig:imex_parareal_512_close} and \ref{fig:imex_parareal_2048_close} we show a grid of convergence plots and stability-convergence overlay plots for IMEX Parareal integrators with $\Nsteps=512$ and $\Nsteps=2048$. Due to the limited stability of the Parareal methods near the origin, we magnify the axes in comparison to our plots for IMEX-RK methods. %
	Three additional figures for Parareal methods with different coarse and fine integrators are shown in Appendix \ref{ap:additional_figures}, and additional plots with further methods and different axes are available in \cite{t_buvoli_2021_4513662}.

\subsubsection{General remarks on stability and convergence} We make several remarks based on the stability plots shown in Figures \ref{fig:imex_parareal_512_close}, \ref{fig:imex_parareal_2048_close}, %
 and the additional ones contained in \cite{t_buvoli_2021_4513662}. These remarks are meant to be general and should hold for parareal methods constructed using a larger range of IMEX-RK families than the four discussed in this work.

\begin{enumerate}
	
	\item There are many regions in the $(z_1, z_2)$ plane where stability and convergence regions do not overlap (in Figures \ref{fig:imex_parareal_512_close} and \ref{fig:imex_parareal_2048_close} these regions appear in a light blue color). Therefore a Parareal method that does not take $K=\Nprocs$ can be unstable even for ($z_1, z_2$) pairs that are inside the convergence region. If this Parareal method is iterated over many blocks, or if the instability is large, then the final output will be unusable. Conversely, lack of convergence paired with stability means that the solution will remain bounded. However the accuracy for these $(z_1, z_2)$ pairs will never surpass that of the coarse integrator unless the number of iterations is large; in many cases as large as $K=\Nprocs$. 
	
	\item Convergence regions for Parareal change drastically depending on the parameters. Overall we see that the choice of coarse integrator and the number of processors have the most prominent effect on the convergence regions near the origin. In particular, the shape of the convergence region is primarily driven by the choice of coarse integrator (See the following remark, and the additional figures in Appendix \ref{ap:additional_figures}), and the size of the region contracts if $\Nsteps$ is increased while $\NFprop$ is kept constant. For a fixed $\Nsteps$, the size of the convergence region near the origin is roughly inversely proportional to $\NFprop$ (or equivalently, proportional to the ratio of the number of processors $\Nprocs$ and the block size $\Nsteps$). Though one may hope to attain the fastest speedup using a coarse solver with large timesteps, these observations suggest that it is important to carefully balance speedup with the size of the solution spectrum.
	\item The choice of fine integrator has almost no effect on the shape or size of the convergence region near the origin. We can understand this by first noticing that the convergence region depends on $F$ through the term $|G-F|$ in (\ref{eq:parareal_iteration_inf_norm}). Now, suppose that we select a coarse integrator of order $\gamma$, and a fine integrator of order $\gamma+\delta$, where $\delta \ge 0$. If we define $\lambda=z_1+z_2$, it follows that
		\begin{align*}
			G = \exp(i\lambda)+ \mathcal{O}\left( |\lambda|^\gamma \right),	 &&
			F(\delta) = \exp(i\lambda)+ \mathcal{O}\left( |\lambda|^{\gamma+\delta} \right).
		\end{align*}
		Therefore $|G-F(\delta)| = |G-F(0)| + \mathcal{O}\left( |\lambda|^\gamma \right)$. This implies that for any choice of fine integrator, the convergence rates will be nearly identical near the $(z_1, z_2)$ origin where $\lambda$  is small.
	
	\item Stability regions are non-trivial and depend on more parameters than convergence regions. The most significant effect on stability regions is due to the choice of coarse integrator, the number of processors, and the number of Parareal iterations.	In the next subsection, we make a few additional comments about the stability of the IMEX integrators tested in this work.

\end{enumerate}

\subsubsection{Remarks regarding IMEX Parareal}

Overall, IMEX Parareal methods are not well-suited for non-diffusive equations. Our main observation is that the stability region splits in two along the $z_1$ axis; see the bottom right diagrams in Figures \ref{fig:imex_parareal_512_close} and \ref{fig:imex_parareal_2048_close} as an example. For all ten IMEX-RK pairings that form the basis of this work, we consistently find the following patterns for the stability regions near the origin:
	\begin{itemize}
		\item Stability regions grow larger if $\Nsteps$ is increased while keeping $\NFprop$ constant.
		\item For fixed $\Nsteps$, stability regions grow if $\NFprop$ is decreased (or equivalently $\Nprocs$ is increased).
		\item  Stability regions do not monotonically increase with iteration count. Instead the regions initially contract and  separate as $K$ increases.	%
	\end{itemize}
Nevertheless, IMEX Parareal methods can still be effective for solving non-diffusive equations under three conditions:
\begin{enumerate}
	\item {\bf Select a stable method pairing.} This is non trivial and many pairings lead to an IMEX Parareal method with no stability along the $z_1$ axis for small $K$. Of the ten possible coarse-fine integrator pairings, only the ones with IMEX-RK1 or IMEX-RK3 as the coarse integrator were stable. Selecting IMEX-RK2 or IMEX-RK4 as the coarse integrator produced an unstable Parareal method for all the parameters we considered. More generally, this suggests that RK pairs should be specially constructed to maximize stability.
	\item {\bf Keep the iteration counts low, choose large block lengths, and avoid using too few processors.} Increasing the number of Parareal iterations causes the stability region to separate along the $z_1$ axis. After looking across a wide range of parameters, we find that the stability regions of the IMEX integrators consistently get worse as the number of iterations initially increases. This is especially true if we consider small block lengths $\Nsteps$ or small $\Nprocs$. This suggests that adaptive implementations of IMEX Parareal on non-diffusive problems will likely lead to an unstable method if the residual tolerance is set too low. However, if one fixes the maximum number of iteration $K$, selects a large block length $\Nsteps$, and uses a sufficient number of processors $\Nprocs$, then the methods can be effective.

		We note that choosing large block lengths with many processors will limit the maximal theoretical speedup of the parareal method. However, for non-diffusive equations, optimizing purely for parallel speedup will lead either to an unstable method, or to a parareal method with a very small convergence region that renders it less efficient compared to the serial fine integrator. Moreover, the important metric is not parallel speedup, but rather efficiency (error vs time) of the parareal integrator compared to its serial fine integrator. To properly determine whether a parareal configuration is efficient, one needs to carefully balance stability, speedup, and convergence relative to the spectrum of the ODE problem that is being solved.
	
	\item {\bf Avoid problems with broad spectrums}. Stable pairings of IMEX integrators possess good stability and convergence near the $(z_1, z_2)$ origin. However for all the RK pairings we tested, the convergence regions do not extend far along the $z_1$ axis. For moderately sized $z_1$ we consistently see a region of good stability that is paired with non-contractive convergence. These observations suggest that IMEX Parareal methods will converge slowly on solutions where the energy is concentrated in these modes. 
 We note that the convergence region is restored as $z_1$ gets sufficiently large (outside the range of the figures contained in the paper). However, in these regions the coarse and fine IMEX integrators both exhibit heavy damping, therefore, rapid convergence to the fine integrator is not a sign of good accuracy; see remark \ref{remark:parareal_dispersive_dissipative} in Section \ref{subsec:convergence}.

\end{enumerate}

\subsection{Accuracy regions for the non-diffusive Dahlquist equation}

To supplement our stability plots we also consider the accuracy regions for IMEX integrators. The accuracy region $\A_\epsilon$ shows the regions in the ($z_1$, $z_2$) plane where the difference between the exact solution and the numerical method is smaller than $\epsilon$. The accuracy region for an IMEX method is typically defined as
	\begin{align}
		\A_\epsilon = \left\{ (z_1, z_2) \in \mathbb{R} : |R(iz_1, iz_2) - \exp(iz_1, iz_2) | \le \epsilon \right\},	
	\end{align}
	where $R(\zeta_2,\zeta_2)$ is the stability function of the method. Since Parareal methods advance the solution multiple timesteps, we scale the accuracy regions by the total number of timesteps in a block so that
	\begin{align}
		\A_\epsilon = \left\{ (z_1, z_2) \in \mathbb{R} : |R(i \Nsteps z_1, i \Nsteps z_2) - \exp(i \Nsteps z_1, i\Nsteps z_2) | \le \epsilon \right\}.
	\end{align}
	Under this scaling, the accuracy region of a fully converged Parareal method with $K = \Nprocs$ is 
		\begin{align}
		\A_\epsilon = \left\{ (z_1, z_2) \in \mathbb{R} : |R_F(i z_1, i z_2)^{\Nsteps} - \exp(i \Nsteps z_1, i\Nsteps z_2) | \le \epsilon \right\}.	
	\end{align}
	where $R_F(\zeta_1, \zeta_2)$ is the stability function of the fine integrator; this is the accuracy region of a method that consists of $\Nsteps$ steps of the fine integrator.

\begin{figure}[h!]	
	\centering
	\includegraphics[width=0.98\linewidth,trim={10cm 20cm 8cm 9cm},clip]{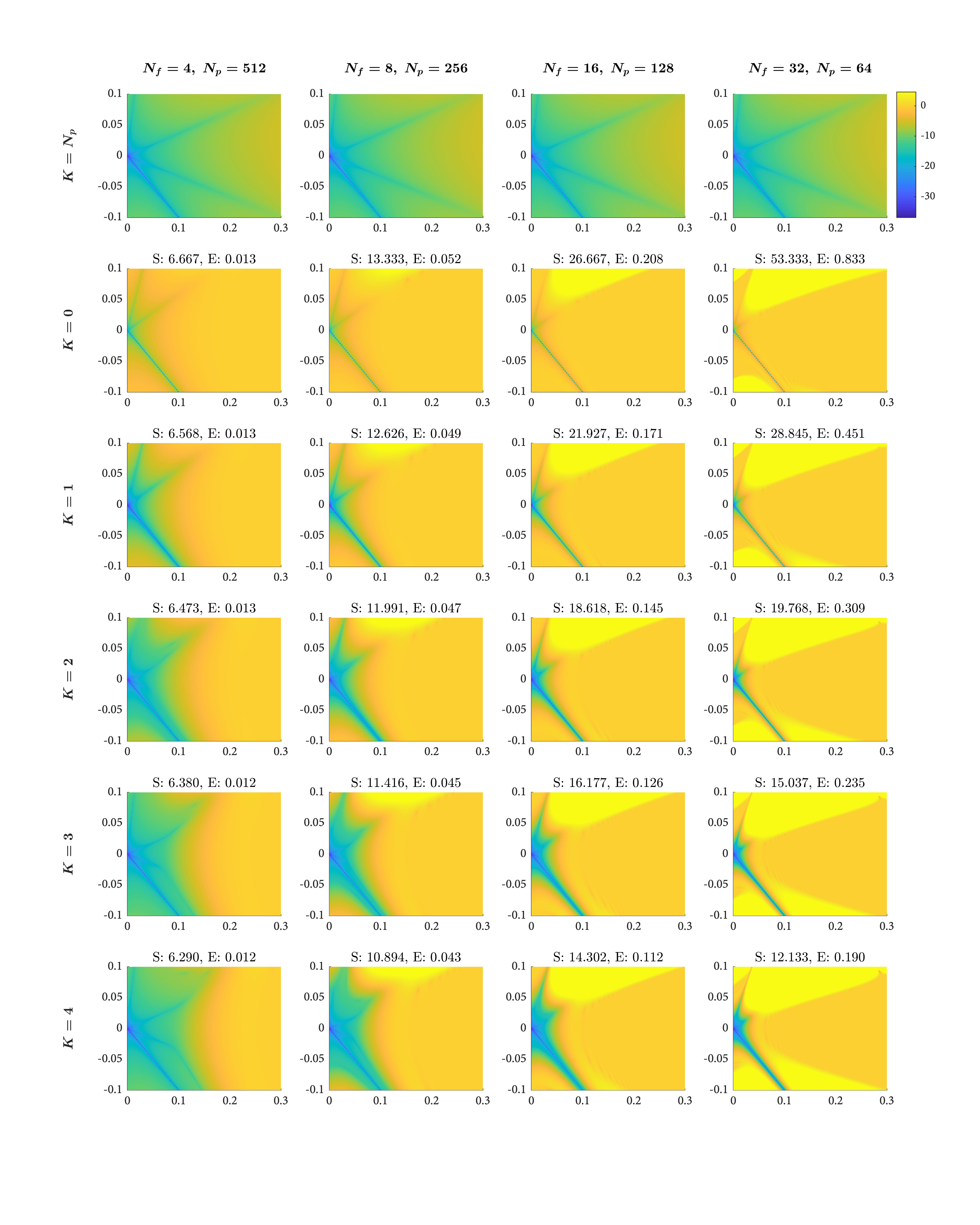}	
	
	\caption{Accuracy regions for Parareal with a block size of $\Nsteps=2048$ and IMEX-RK3, IMEX-RK4 as the coarse and fine integrator. These plots complement the stability plots in Figure \ref{fig:imex_parareal_2048_close}. The top row shows the accuracy of the fine integrator, the second row shows the coarse integrator, and all subsequent rows show parareal methods with an increasing number of iterations. Color represents the log of the absolute value of the error. }
	\label{fig:accuracy_imex_parareal_2048_close}	
\end{figure}

	In Figure \ref{fig:accuracy_imex_parareal_2048_close} we show the accuracy regions for IMEX Parareal methods with $\Nsteps=2048$. The diagrams highlight the importance of the different regions in our stability plots. First, we can clearly see fast convergence to the fine solution inside contractive regions where $\cRate \ll 1$. Second, for any ($z_1$, $z_2$) pairs that are inside the stability region but outside the convergence region, we see that increasing the number of Parareal iterations does not improve the accuracy of the solution beyond that of the coarse integrator. Finally, in the regions with large instabilities the solution is no longer useful. However many of the instabilities that lie inside the convergence regions are so small in magnitude that they do not appear prominently in the accuracy plots, since they will only affect accuracy after many steps of Parareal.

\section{Numerical experiments}\label{sec:numerical_experiments}

The aim of our numerical experiments is to validate the results from linear stability analysis and demonstrate that we can obtain meaningful parallel speedup with IMEX-RK Parareal on a nonlinear dispersive equation. In particular, we show that: 
\begin{enumerate}
    \item Taking too many Parareal iterations leads to an unstable method; however, selecting larger block sizes $\Nsteps$ for a fixed $\NFprop$ increases the range of stable choices for the iteration count $K$. 
    \item Decreasing the number of processors for a fixed block size $\Nsteps$ and fixed $\NGprop$ eventually leads to an unstable method, even for timesteps where the coarse integrator is stable.
    \item For the parameters considered in this work, Parareal methods with adaptive iteration count $K$ are effective so long as one limits the maximum number of iterations.
\end{enumerate}

For our model nonlinear problem we select the one dimensional nonlinear Schr\"{o}dinger (NLS) equation
\begin{align}
		\begin{aligned}
			& iu_t + u_{xx} + 2|u|^2u = 0, \\
			& u(x, t=0) = 1 + \tfrac{1}{100} \exp(ix / 4), \quad x \in [-4\pi, 4\pi].
			\label{eq:nls_problem}
		\end{aligned}
	\end{align}
We equip NLS with periodic boundary conditions and use the method of lines with a fixed spatial grid. 	For the spatial component, we use a Fourier spectral discretization with 1024 points. In the time dimension we apply IMEX Parareal methods that treat the linear spatial derivative term implicitly and the nonlinear term explicitly. We solve the equation in Fourier space where the discrete linear derivative operators are diagonal matrices and the implicit solves amount to multiplications with a diagonal matrix. %
	
	The solution of the problem has a subtle behavior when used as a test for a parallel-in-time method. The initial condition is a perturbation of the plane wave solution $u(t,x) = \exp(-2it)$. For small time, the solution is smooth and the perturbation experiences exponential growth that is well-described by the PDE's local linearization around the plane wave solution \cite{buvoli2013rogue}. Then, around  $t=10$, the nonlinear terms become dominant and cause spectral broadening. 
	
	We integrate (\ref{eq:nls_problem}) out to time $t=15$ using $\Nstepstot = 2^p$ timesteps where ${p=7, \ldots, 18}$. For a Parareal method with block size $\Nsteps$ we only show data points for $\Nstepstot \ge \Nsteps$ since one cannot run Parareal with fewer than $\Nsteps$ timesteps. Note that when $\Nstepstot > \Nsteps$, this amounts to computing the final solution using multiple parareal blocks where $\Nblocks = \Nstepstot / \Nsteps$. For brevity we only consider Parareal integrators with IMEX-RK3 as the coarse integrator and IMEX-RK4 as the fine integrator and always take the number of coarse steps $\NGprop = 1$. In all our numerical experiments, we also include a serial implementation of the fine integrator.
	
	Finally, for all the plots shown in this section, the relative error is defined as $\|\mathbf{y}_{\text{ref}} - \mathbf{y}_{\text{method}} \|_\infty / \| \mathbf{y}_{\text{ref}} \|_\infty$ where $\mathbf{y}_{\text{ref}}$ is a vector containing the reference solution in physical space and $\mathbf{y}_{\text{method}}$ is a vector containing the output of a method in physical space. The reference solution was computed by running the fine integrator (IMEX-RK4) with $2^{19}$ timesteps, and the relative error is always computed at the final time $t=15$.

\subsection{Varying the block size $\Nsteps$ for fixed $\NFprop$ and $\NGprop$ } %
In our first numerical experiment we show that increasing the block size $\Nsteps$ for fixed $\NFprop$ and $\NGprop$ allows for Parareal configurations that remain stable for an increased number of iterations $K$.  Since we are fixing $\NFprop$ and $\NGprop$, we are increasing $\Nprocs$ to obtain larger block sizes. Therefore, this experiment simultaneously validates the improvement in stability seen when comparing the $i$th column of Figure \ref{fig:imex_parareal_2048_close} with the $i$th column of Figure \ref{fig:imex_parareal_512_close}, along with the decrease in stability seen when moving downwards along any column of Figure \ref{fig:imex_parareal_512_close} or \ref{fig:imex_parareal_2048_close}. 

\begin{figure}[h!]
	
	\begin{center}
		\includegraphics[height=0.4\linewidth,trim={0 0 0 0},clip]{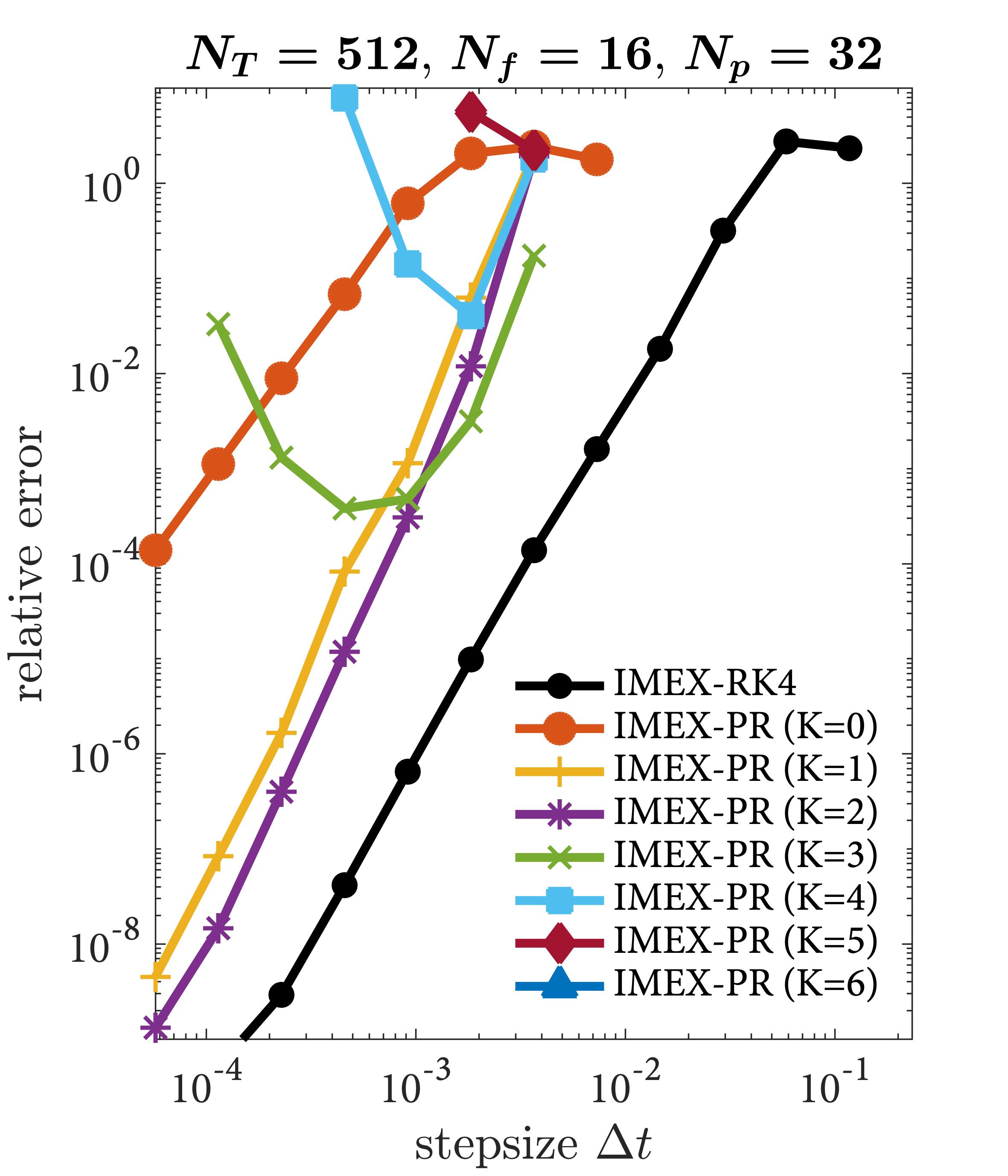}
		\includegraphics[height=0.4\linewidth,trim={6cm 0 0 0},clip]{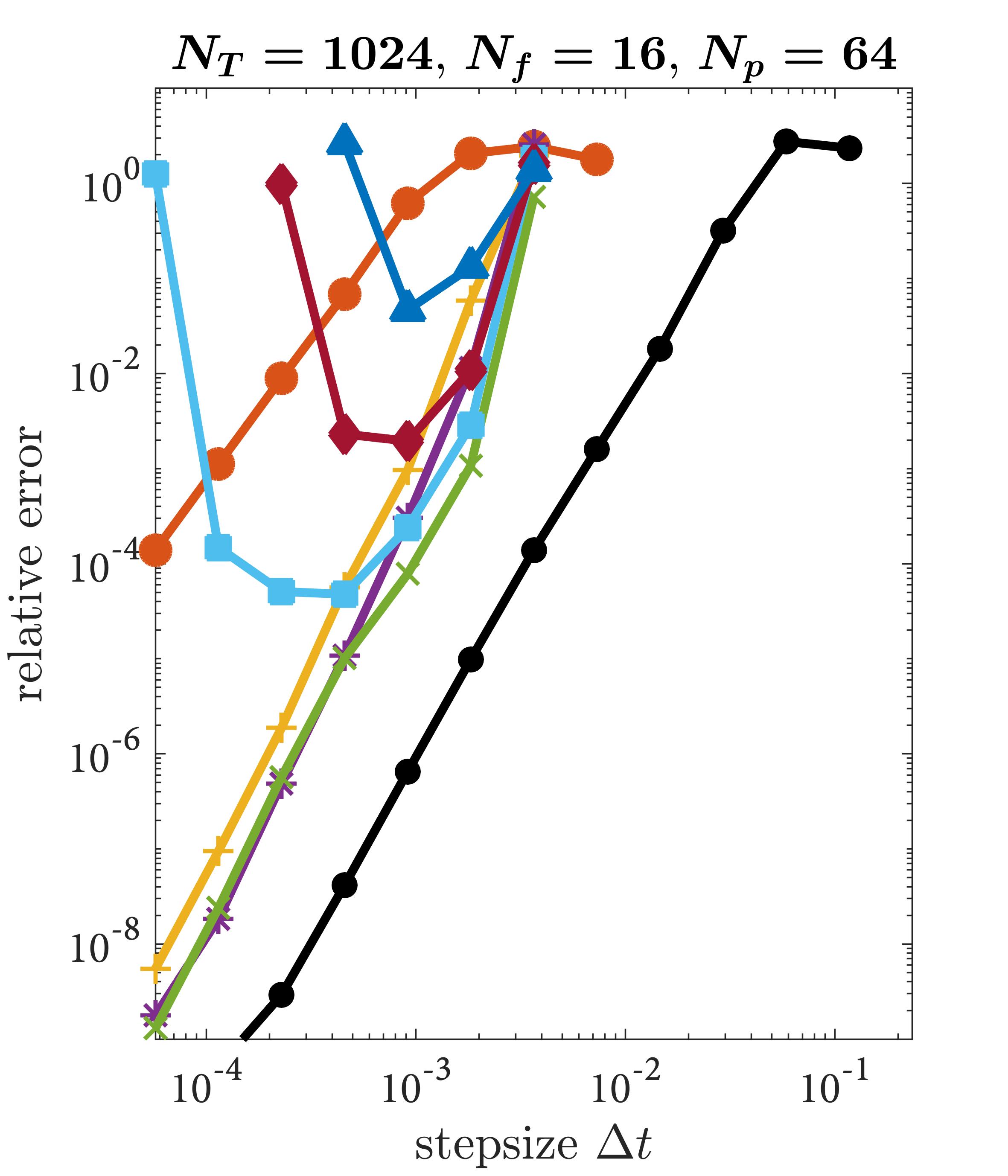}
		\includegraphics[height=0.4\linewidth,trim={6cm 0 0 0},clip]{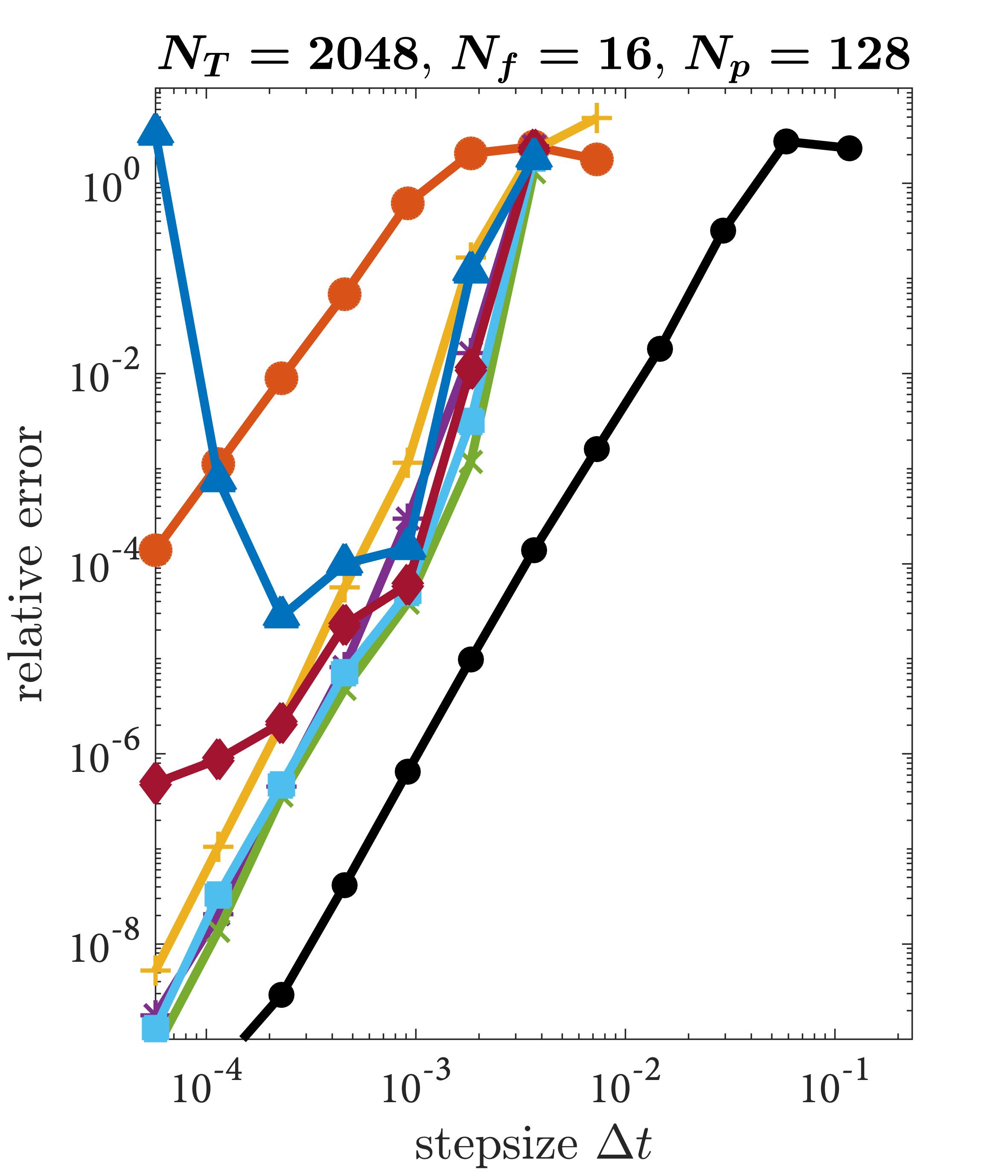}
	\end{center}

	\caption{Accuracy versus stepsize plots for the Parareal method with IMEX-RK3, IMEX-RK4 as the coarse and fine integrator. The block size $\Nsteps$ for the left, middle, and right plots is respectively 512, 1024, and 2048.  The black line shows the serial fine integrator, while the colored lines represent Parareal methods with different values of iteration count $K$. Note that the Parareal configuration with $\Nsteps = 512$ and $K=6$ did not converge for any of the timesteps. }
	\label{fig:nls_error_vs_stepsize}
\end{figure}

In Figure \ref{fig:nls_error_vs_stepsize} we present plots of relative error versus stepsize for three Parareal configurations with block sizes $\Nsteps = 512$, $1024$, or $2048$. Each of the configurations takes a fixed number of Parareal iterations $K$ and has ${\NGprop=1}$, ${\NFprop=16}$. The stability regions for the Parareal configurations with $\Nsteps=512$ and ${\Nsteps=2048}$ are shown in the third columns of Figures  \ref{fig:imex_parareal_512_close} and \ref{fig:imex_parareal_2048_close}. From linear stability we expect that the two Parareal methods will respectively become unstable if $K > 3$ and  $K > 4$. The experiments with $\Nsteps=512$ align perfectly with the linear stability regions. For the larger block size, the instabilities are milder and we would need to take many more Parareal blocks for them to fully manifest. Nevertheless the methods fail to become more accurate for $K=4$ and start to diverge for  $K > 4$. Overall, the results confirm that increasing the block size by increasing $\Nprocs$ leads to an improvement in stability that allows for a larger number of total Parareal iterations. 

Finally, we note that an alternative strategy for increasing $\Nsteps$  is to increase $\NFprop$ while keeping $\Nprocs$ constant. However, we do not consider this scenario since increasing $\NFprop$ will lead to a method with a significantly smaller convergence region and no better stability (for example, compare column 2 of Figure \ref{fig:imex_parareal_512_close} with column 4 of Figure \ref{fig:imex_parareal_2048_close}).

\subsection{Varying the number of processors $\Nprocs$ for a fixed $\Nsteps$ and $\NGprop$}

\begin{wrapfigure}{r}{0.4 \linewidth}
  	\vspace{-1.5em}
  	\centering
	\includegraphics[width=\linewidth,trim={1cm 1cm 5cm 3cm},clip]{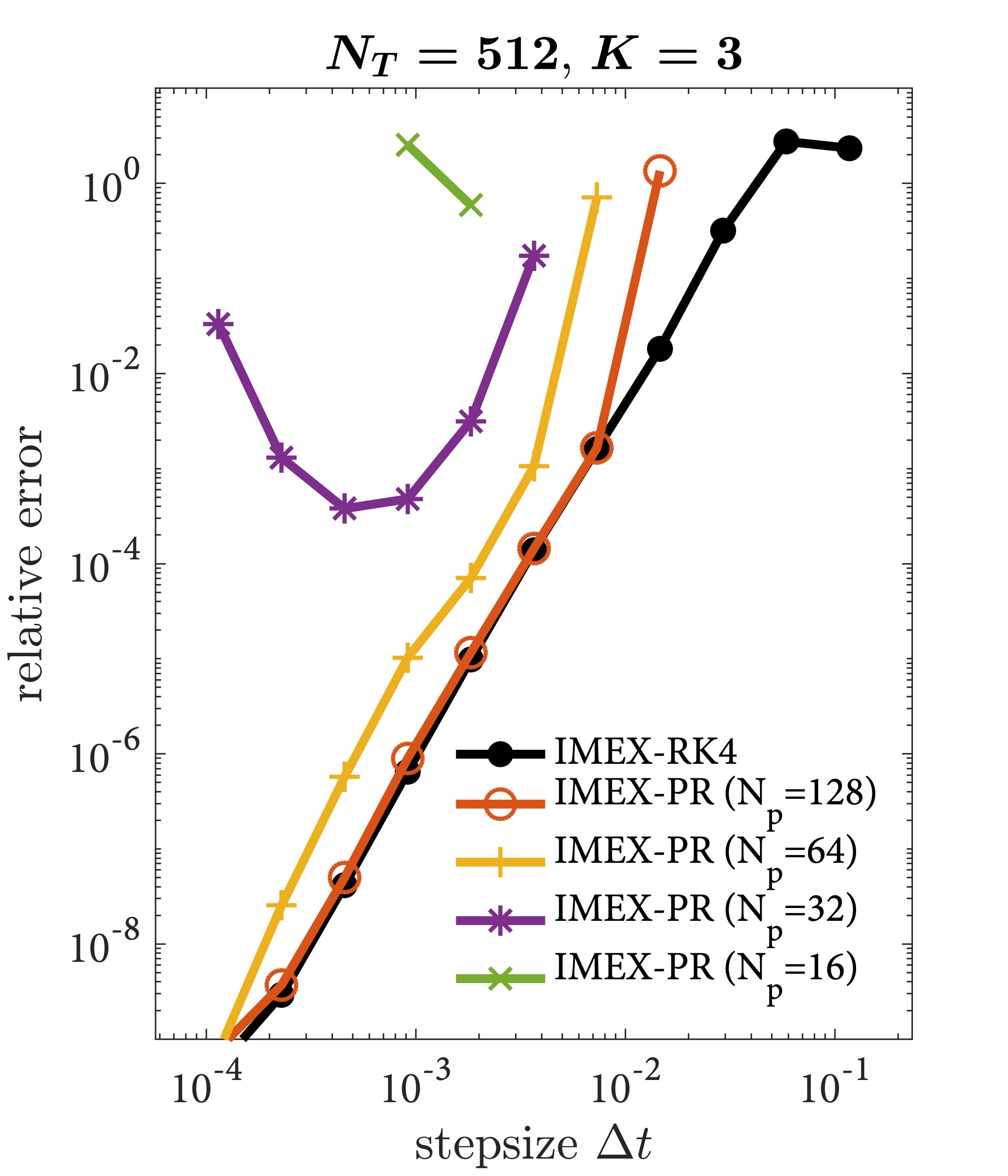}
	\caption{Variable $\Nprocs$ results.}
	\vspace{-1.5em}
	\label{fig:nls_variable_NFsteps}
\end{wrapfigure}

In our second numerical experiment we show that decreasing the number of processors $\Nprocs$ (or equivalently increasing the the number of fine steps $\NFprop$) for a fixed $\Nsteps$ and $\NGprop$ will lead to an unstable Parareal method. This experiment validates the stability changes that occur along any row of Figure \ref{fig:imex_parareal_512_close} or Figure \ref{fig:imex_parareal_2048_close}.

 For brevity we only consider four Parareal methods with $\NGprop = 1$, $\Nsteps = 512$, $K=3$, and $\Nprocs = 16$, $32$, $64$, or $128$. The linear stability regions for each of these methods are shown in the fourth row of Figure \ref{fig:imex_parareal_512_close}. Only the Parareal method with $\Nprocs=128$ is stable along the entire $z_2$ axis. The method with $\Nprocs=64$ has a mild instability located inside its convergence region, and methods with $\Nprocs=32$ or $\Nprocs=16$ have large instabilities.

In Figure \ref{fig:nls_variable_NFsteps}, we show an accuracy versus stepsize plot that compares each of the four Parareal configurations (shown in colored lines) to the serial fine integrator (shown in black). The total number of steps $\Nstepstot$ is given by $2^p$ where for Parareal ${p=9, \ldots, 18}$, while for the serial integrator ${p=7,\ldots, 18}$. We can clearly see that decreasing $\Nprocs$ leads to instability even at small stepsizes. Note that for $\Nprocs=64$ the instability is so small that it does not affect convergence in any meaningful way. 

It is important to remark that the largest stable stepsize for any Parareal method is restricted by the stability of its coarse integrator. Since we are taking ${\NGprop = 1}$, the number of coarse timesteps per block is $\Nprocs$. Therefore, a Parareal configuration with a smaller $\Nprocs$ takes larger coarse stepsizes and requires a smaller $\Delta t$ to remain stable. This effect can be see in Figure \ref{fig:nls_variable_NFsteps}, since it causes the rightmost point of the parareal convergence curves to be located at smaller stepsizes $\Delta t$ for methods with smaller $\Nprocs$. What is more interesting however, is that the Parareal configurations with fewer processors are unstable, even when the stepsize is lowered to compensate for the accuracy of the coarse solver. In other words, instabilities form when the difference in accuracy between the coarse and fine solver is too large, even if the coarse solver is sufficiently stable and accurate on its own.

\subsection{Efficiency and adaptive $K$}

In our final numerical experiment we first compare the theoretical efficiency of several IMEX Parareal configurations and then conduct a parallel experiment using the most efficient parameters. We conduct our theoretical efficiency analysis using Parareal methods with $\Nsteps=2048$, $\NFprop=16$, $\NGprop=1$, and $K=1, \ldots, 6$. As shown in the third column of Figure \ref{fig:accuracy_imex_parareal_2048_close}, these configurations possess good speedup and stability regions when $K \le 3$. To determine the {\em theoretical runtime} for Parareal (i.e the runtime in the absence of any communication overhead), we divide the runtime of the fine integrator by the Parareal speedup that is computed using (\ref{eq:parareal_speedup2}).

\begin{figure}
	
	\begin{center}
		\includegraphics[width=0.47\linewidth]{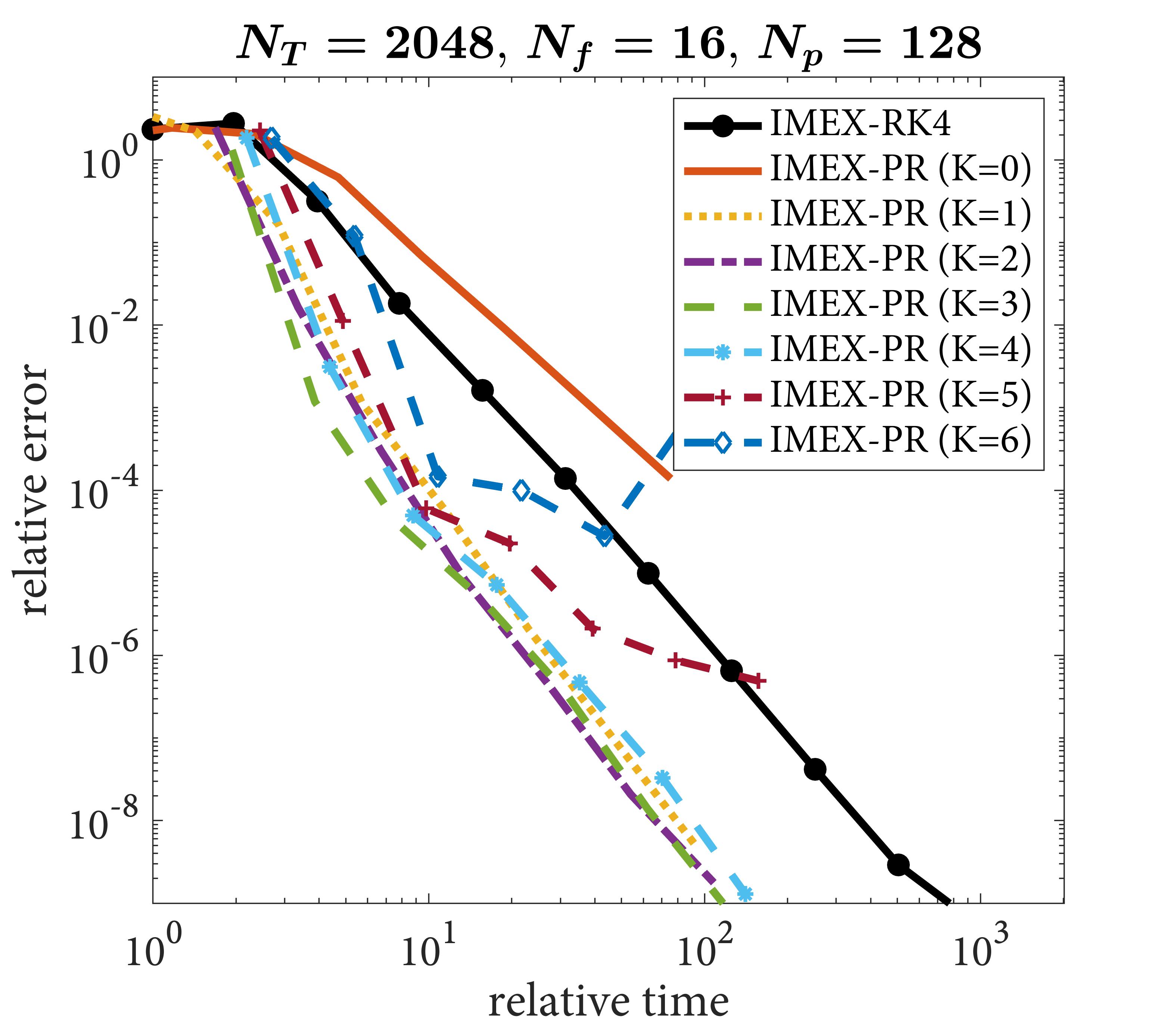}
		\includegraphics[width=0.47\linewidth]{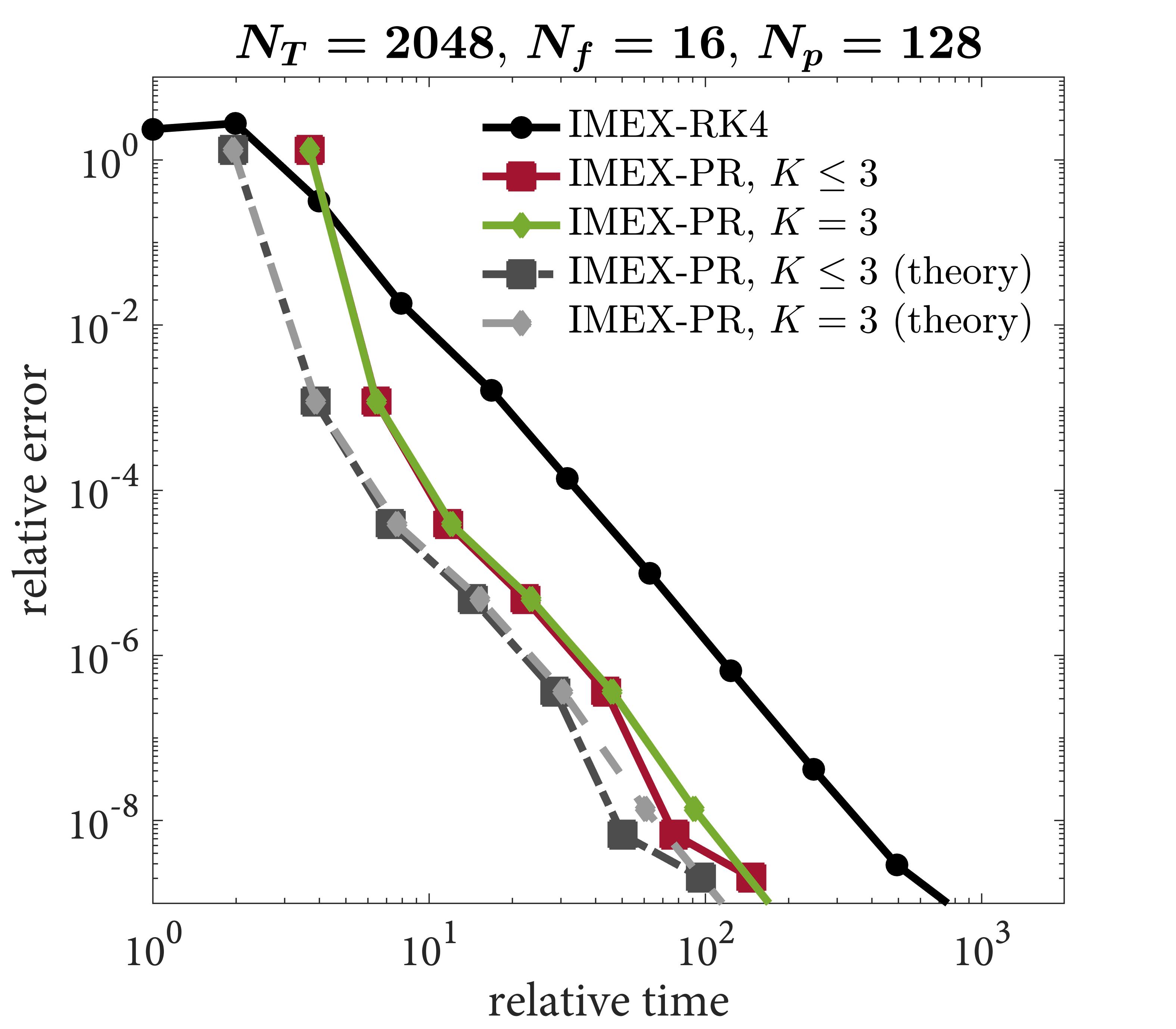}

		\begin{minipage}[t]{0.47\linewidth}	
			\centering (a)
		\end{minipage}
		\begin{minipage}[t]{0.47\linewidth}	
			\centering (b)
		\end{minipage}
	\end{center}

	\caption{Relative error versus computational time for the NLS equation solved using an IMEX Parareal methods with $\Nsteps=2048$ and $\Nprocs=128$. The left plot (a) compares the theoretical running times of seven Parareal methods that each take a different number of iterations $K$ per block. The right plot (b) compares the real-world running times of the Parareal method with $K=3$, and a Parareal method with adaptive controller where $K\le 3$. We also show the theoretical running times of the two methods in grey to highlight the losses due to communication. In both plots, the black line shows the fine integrator that is run in serial. All times have been scaled relative to the fine integrator at the coarsest timestep.}
	\label{fig:nls_accuracy_vs_stepsize}
\end{figure}

In Figure \ref{fig:nls_accuracy_vs_stepsize}(a) we show plots of relative error versus theoretical runtime for the seven parareal configurations. The total number of steps $\Nstepstot$ is given by $2^p$ where for Parareal ${p=12, \ldots, 18}$, while for the serial integrator ${p=7,\ldots, 18}$. Note that the running times for the fine integrator and the parareal configuration with $K=0$ (i.e. the serial coarse method) measure real-world efficiency since there are no parallelization possibilities for these methods. 

The efficiency plots demonstrate that it is theoretically possible to achieve meaningful parallel speedup using IMEX Parareal on the nonlinear Schr\"{o}dinger equation. Moreover, amongst the seven Parareal configurations, the one with $K=3$ is the most efficient over the largest range of timesteps. However, communication costs on real hardware are never negligible, and real-world efficiency will depend heavily on the underlying hardware and the ODE problem. To validate the practical effectiveness of IMEX Parareal, we ran the most efficient configuration with $K=3$ in parallel on a distributed memory system with 128 processors \footnote{The numerical experiments were performed on the Cray XC40 "Cori" at the National Energy Research Scientific Computing Center using four 32-core Intel "Haswell" processor nodes.  The Parareal method is implemented as part of the open source package LibPFASST available at \url{https://github.com/libpfasst/LibPFASST}.}. We also tested in identical Parareal configuration with an adaptive controller for $K$ that iterates until either $K \ge K_{max}$, or a residual tolerance of $1\times 10^{-9}$ is satisfied. Unsurprisingly, it was necessary to restrict the maximum number of adaptive Parareal iterations to $K_{max}=3$ or the adaptive controller caused the method to become unstable.

In Figure \ref{fig:nls_accuracy_vs_stepsize}(b) we show plots of relative error versus parallel runtime, and in Table \ref{tab:nls-speedup} we also include the corresponding speedup for the two Parareal methods. Even on this simple one dimensional problem, we were able to achieve approximately a ten-fold real-world speedup relative to the serial IMEX-RK4 integrator. This is very encouraging since the ratio between the communication and timestep costs is larger for a one dimensional problem. Our results also show that there is not much noticeable difference between the parareal method with fixed $K$ and the method with adaptive $K$, except at the finest timesteps where the adaptive implementation is able to take fewer iterations. %

\begin{table}
	\centering
	\begin{tabular}{l|llll}
		$\Nstepstot$ & AS ($K=3$)~ & TS ($K=3$)~ & AS ($K\le 3$)~ & TS ($K\le 3$)\\ \hline
		4096 & 8.55 &		16.18	&	8.54		&	16.18  	\\
	8192 & 		9.75 &		16.18	&	9.74 	& 	16.18	\\
	16384 & 	10.25 &		16.18	&	10.55 	& 	17.01	\\
	32768 & 	10.56 &		16.18	&	11.06	& 	17.16	\\ 
	65536 & 	10.74 &		16.18	&	11.27 	& 	17.24	\\
	131072 & 	10.76 &		16.18	&	12.69 	&	19.54	\\
	262144 & 	10.98 &		16.18	&	13.57	&	20.61
	\end{tabular}
	\vspace{1em}
	\caption{ Achieved speedup (AS) and theoretical speedup (TS) for the two Parareal configurations shown in Figure \ref{fig:nls_accuracy_vs_stepsize}(b).}
	\label{tab:nls-speedup}
\end{table}

\section{Summary and conclusions} 
\label{sec:summary}

We have introduced a methodology for categorizing the convergence and stability properties of a Parareal method with pre-specified parameters. By recasting the Parareal algorithm as a one-step RK method with many stages, we are able to combine classical stability analysis with a simple bound on the norm of the Parareal iteration matrix. The resulting stability convergence overlay plots highlight the key characteristics of a Parareal method including: regions of fast and slow convergence, stable regions where convergence does not occur, and regions where instabilities will eventually contaminate the method output.

By searching through a wide range of IMEX Parareal methods, we were able to identify several stable  configurations that can be used to solve dispersive equations. Moreover, each of the configurations possessed the same characteristics: low iteration counts $K$, large block sizes $\Nsteps$, and a large number of processors $\Nprocs$. We also observed that the coarse integrator is the most important factor that determines whether a Parareal method is stable, and a bad choice can single-handedly lead to an unstable method regardless of the other parameters.

More broadly, we see that convergence and stability regions are highly nontrivial and depend heavily on the parameters. It is clear that one cannot arbitrarily combine coarse and fine integrators and expect to obtain a good Parareal method for solving dispersive equations. The same lesson also applies to all Parareal parameters since serious instabilities can form by arbitrarily changing the number of iterations, the block size, or the number of processors. 

Finally, we remark that the analysis presented in this work can be reused to study the properties of any Runge-Kutta Parareal method on the more general partitioned Dahlquist problem that represents both dispersive and diffusive equations. However, many of the conclusions and properties that we found are specific to IMEX methods and will not hold for different method families or for different problem types.
 
\section*{Acknowledgements}
The work of Buvoli was funded by the National Science Foundation, Computational Mathematics Program DMS-2012875.

The work of Minion was supported by the U.S. Department of Energy, Office of Science,
Office of Advanced Scientific Computing Research, Applied Mathematics program
under contract number DE-AC02005CH11231. Part of the simulations were performed
using resources of the National Energy Research Scientific Computing Center
(NERSC), a DOE Office of Science User Facility supported by the Office of
Science of the U.S. Department of Energy under Contract No. DE-AC02-05CH11231.

\bibliographystyle{siam}
\bibliography{references_nourl,pint,imex,gen}

\clearpage{}%
\appendix

\section{Infinity norm of the Parareal iteration matrix $\mathbf{E}$}
\label{ap:inf_norm_parareal_iteration_matrix}

Let $\mathbf{A}(\gamma)$ be the lower bidiagonal matrix
\begin{align*}
\mathbf{A}(\gamma) = 
\left[
	\begin{array}{cccc}
		1 & & & \\
		\gamma & 1 & & \\
		& \ddots & \ddots & \\
		& & \gamma & 1
	\end{array}.
\right]
\end{align*}

\begin{lemma}
The inverse of $\mathbf{A}(\gamma)$ is given by
\begin{align*}
	\mathbf{A}^{-1}_{i,j}(\gamma) = 
	\begin{cases}
		(-\gamma)^{i-j} & j \le i \\
		0 & \text{otherwise}
	\end{cases} 
\end{align*}
\begin{proof}
	For convenience we temporarily drop the $\gamma$ so that $\mathbf{A}=\mathbf{A}(\gamma)$, then
	\begin{align*}
			\left( \mathbf{A} \mathbf{A}^{-1} \right)_{ij} = \sum_{k=1}^{\Nprocs+1} \mathbf{A}_{ik}\mathbf{A}_{kj}^{-1} &=  
			\begin{cases}
				0 & j > i\\
				\mathbf{A}_{ii}\mathbf{A}_{ii}^{-1} & i=j \\
				\mathbf{A}_{ii}\mathbf{A}_{ij}^{-1} + A_{i, i-1} A_{i-1,j}^{-1} & j<i \\	
			\end{cases} \\
			&=  
			\begin{cases}
				0 & j > i\\
				1 & i=j \\
				(-\gamma)^{i-j} +  \gamma (-\gamma)^{i-1-j} & j<i	
			\end{cases} \\
			&=  
			\begin{cases}
				1 & i=j, \\
				0 & \text{otherwise}	.
			\end{cases}
	\end{align*}
\end{proof}	
\end{lemma}

\begin{lemma}
\label{lem:matrix-inverse}
The product of $\mathbf{A}(\omega)\mathbf{A}^{-1}(\gamma)$ is
\begin{align*}
	\left( \mathbf{A}(\omega) \mathbf{A}^{-1}(\gamma) \right)_{ij} = 
	\begin{cases}
		0 & j > i \\
		1 & i = j \\
		(-\gamma)^{i-j-1}(\omega - \gamma) & j < i
	\end{cases}
\end{align*}
\begin{proof}
	\begin{align*}
			\left( \mathbf{A}(\omega) \mathbf{A}^{-1}(\gamma) \right)_{ij} 
				& = \sum_{k=1}^{\Nprocs+1} \mathbf{A}_{ik}(\omega)\mathbf{A}_{kj}^{-1}(\gamma) \\
				&=  
			\begin{cases}
				0 & j > i\\
				\mathbf{A}_{ii}(\omega) \mathbf{A}_{ii}^{-1}(\gamma) & i=j \\
				\mathbf{A}_{ii}(\omega) \mathbf{A}_{ij}^{-1} + \mathbf{A}_{i, i-1}(\omega) \mathbf{A}_{i-1,j}^{-1} (\gamma) & j<i
			\end{cases} \\
			&=  
			\begin{cases}
				0 & j > i\\
				1 & i=j \\
				(-\gamma)^{i-j} +  \omega (-\gamma)^{i-1-j} & j<i	
			\end{cases}
	\end{align*}
\end{proof}	
\end{lemma}

\begin{lemma}
The infinity norm of the matrix $M(\omega, \gamma) = \mathbf{I} - \mathbf{A}(\omega)\mathbf{A}^{-1}(\gamma) \in \mathbb{R}^{\Nprocs+1, \Nprocs+1}$ is
\begin{align*}
	\| \mathbf{M}(\omega, \gamma) \|_{\infty} = \frac{1 - |\gamma|^{\Nprocs}}{1 - |\gamma|} |\gamma - \omega|.
\end{align*}
\begin{proof}
	Using Lemma \ref{lem:matrix-inverse}, the $j$th absolute column sum of $M(\omega, \gamma)$ is
	\begin{align*}
		c_j = \sum_{k=j+1}^{\Nprocs+1}	 |(-\gamma)^{k-j-1}(\omega - \gamma)| = \sum_{k=1}^{\Nprocs - j}	 |(-\gamma)^{k}||(\omega - \gamma)|
	\end{align*}
	It follows that $\max_j c_j = c_1$, which can be rewritten as
	\begin{align*}
		\frac{1 - |\gamma|^{\Nprocs}}{1 - |\gamma|} |\gamma - \omega|.
	\end{align*}
\end{proof}
\end{lemma}

\section {Additional Stability and Convergence Overlay Plots}
\label{ap:additional_figures}

The following three figures show stability and convergence overlay plots for Parareal configurations with: $\Nsteps = 2048$, IMEX-RK4 as the fine integrator, and three different coarse integrators. These additional figures supplement Figure \ref{fig:imex_parareal_2048_close} and show the effects of changing the course integrator.

\begin{figure}[h!]
	\centering	
	\includegraphics[width=0.98\linewidth,trim={10cm 20cm 8cm 9cm},clip]{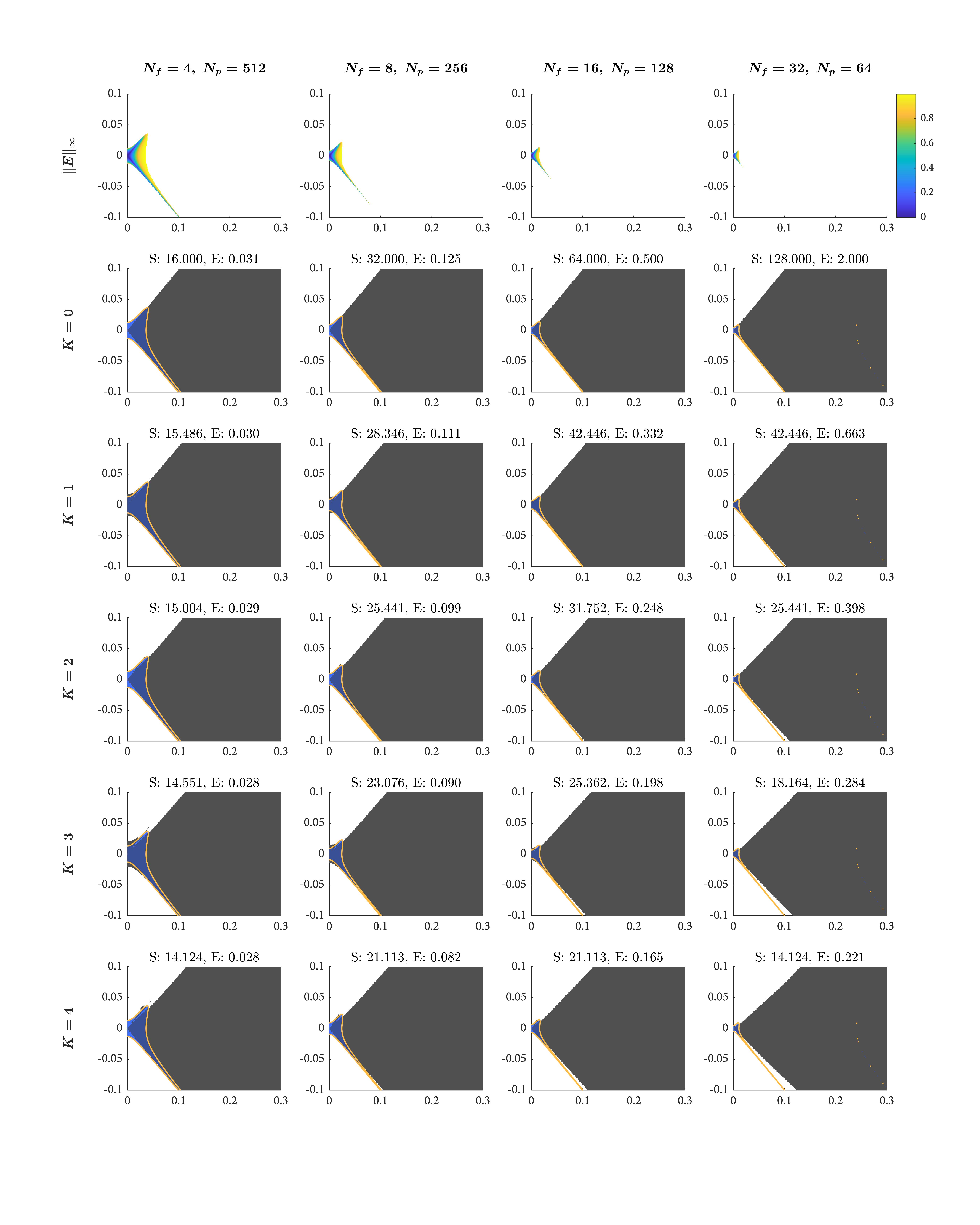}
	
	\caption{Stability and convergence overlay plots for Parareal configurations with a block size of $\Nsteps=2048$ and IMEX-RK1, IMEX-RK4 as the coarse and fine integrator.}
	\label{fig:imex_parareal_2048_close_RK1_RK3}	
\end{figure}

\begin{figure}[h!]
	\centering	
	\includegraphics[width=0.98\linewidth,trim={10cm 20cm 8cm 9cm},clip]{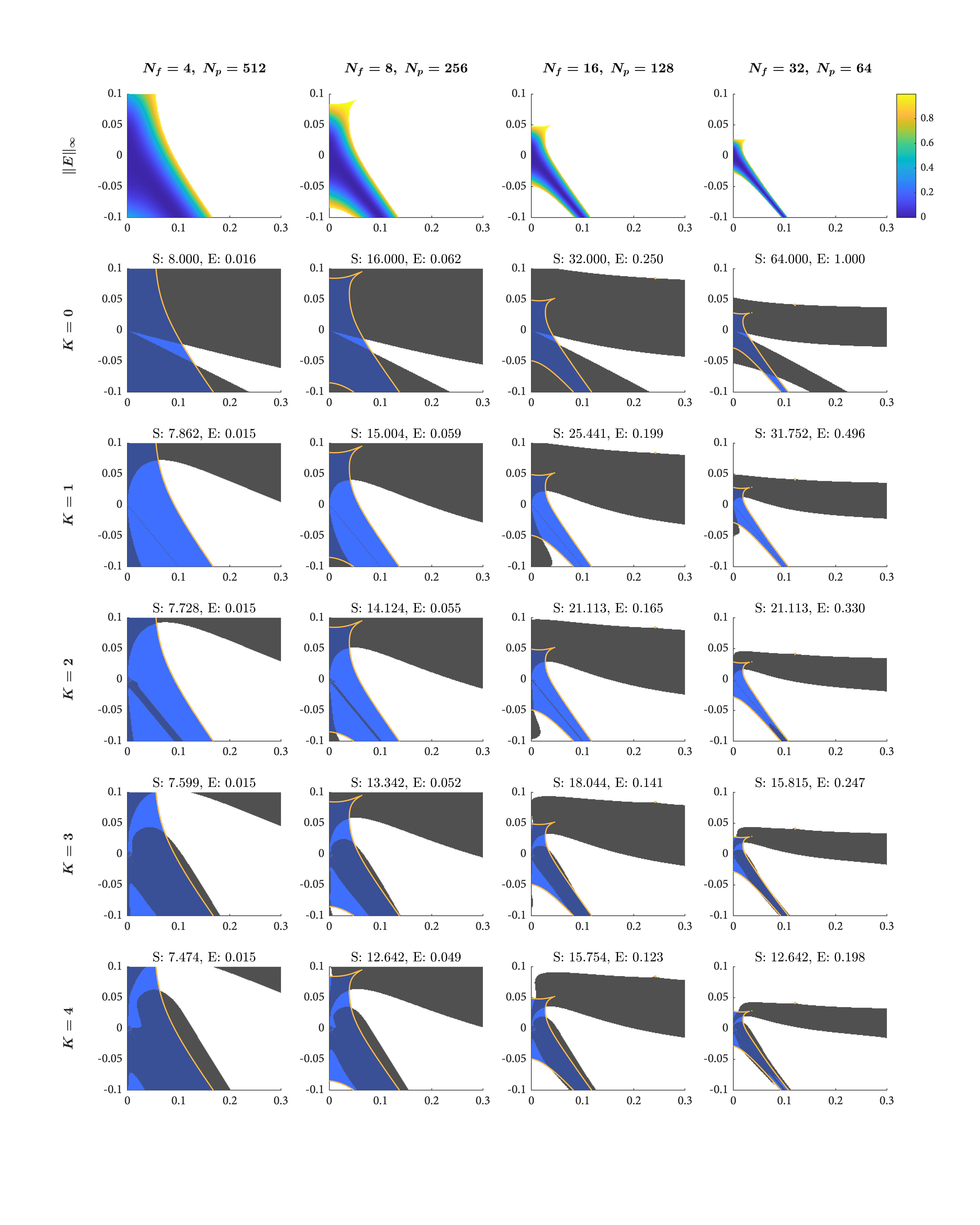}
	
	\caption{Stability and convergence overlay plots for Parareal configurations with a block size of $\Nsteps=2048$ and IMEX-RK2, IMEX-RK4 as the coarse and fine integrator.}
	\label{fig:imex_parareal_2048_close_RK2_RK4}	
\end{figure}

\begin{figure}[h!]
	\centering	
	\includegraphics[width=0.98\linewidth,trim={10cm 20cm 8cm 9cm},clip]{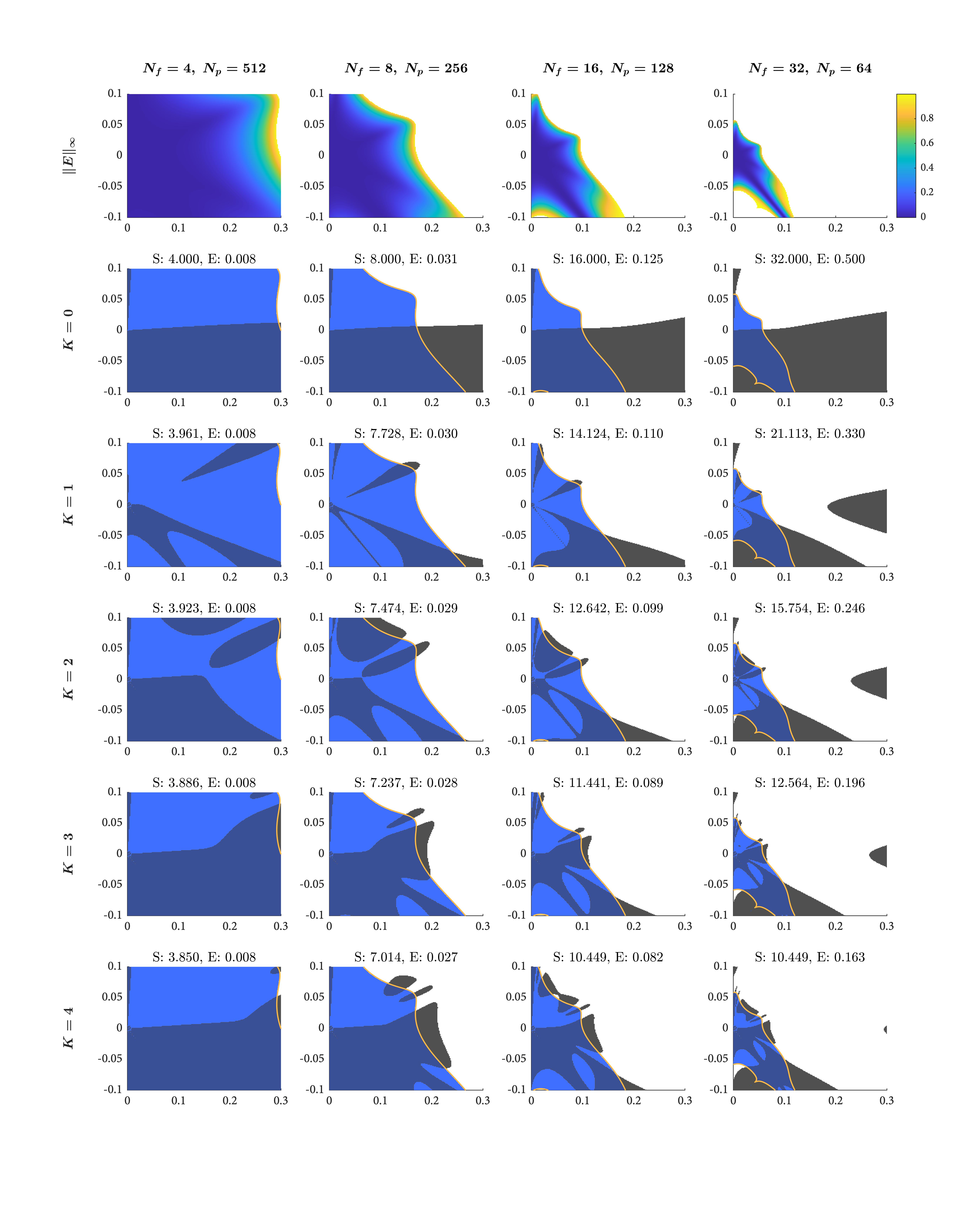}
	
	\caption{Stability and convergence overlay plots for Parareal configurations with a block size of $\Nsteps=2048$ and IMEX-RK4, IMEX-RK4 as the coarse and fine integrator.}
	\label{fig:imex_parareal_2048_close_RK4_RK4}	
\end{figure}

\clearpage{}%

\end{document}